\documentclass[10pt]{article}
\usepackage{geometry,amsmath,amssymb,amsthm,graphicx,epsfig}
\usepackage{subcaption,enumerate,dcolumn,latexsym,epstopdf,color}
\usepackage[graph]{xy}
\usepackage[usenames,dvipsnames,svgnames,table]{xcolor}
\usepackage{ulem}
\numberwithin{equation}{section}
\usepackage{multirow}
\usepackage{tabularx}
\usepackage{mathtools}
\mathtoolsset{showonlyrefs}
\usepackage{pdfsync,comment,a4wide}
\usepackage{enumitem}
\usepackage[titletoc,title]{appendix}
\usepackage{placeins}
\newtheorem{theorem}{Theorem}

\newtheorem{proposition}{Proposition}
\newtheorem{lemma}{Lemma}
\theoremstyle{definition}
\newtheorem{remark}{Remark}
\newtheorem{definition}{Definition}

\newtheorem{assumption}{Assumption}

\renewcommand{\em}[1]{\normalem \em{#1}}
\renewcommand{\emph}[1]{\normalem \emph{#1}}

\newcommand{\nn}{\mathbb{N}} 
\newcommand{\real}{\mathbb{R}} 
 
\newcommand{\norm}[1]{\left\Vert {#1} \right\Vert} 
 
\newcommand{\erl}{\left(-\infty , +\infty\right]} 
\newcommand{\dom}[1]{\mathrm{dom}\,{#1}} 
\DeclareMathOperator*{\argmin}{\arg\!\min}

\newcommand{\dist}{\mathrm{dist}} 
\newcommand{\prox}{\mathrm{prox}} 
 
\newcommand{\act}[1]{\left\langle {#1} \right\rangle} 
\newcommand{\pos}[1]{\left[{#1}\right]_{+}}

\newcommand{\Sequp}[2]{\left\{{#1}^{{#2}}\right\}_{{#2} \in \mathbb{N}}}

\newcommand{\bx}{{\bf x}}
\newcommand{\by}{{\bf y}}
\newcommand{\bw}{{\bf w}}
\newcommand{\bz}{{\bf z}}
\newcommand{\bv}{{\bf v}}
\newcommand{\bu}{{\bf u}}
\newcommand{\bo}{{\bf 0}}
\newcommand{\GGG}{\mathcal{G}}

\usepackage{datetime}
\newdateformat{monthyear}{\monthname[\THEMONTH], \THEYEAR}

\title{A First-Order Method for Nonconvex \\ Simple Bi-Level Optimization}

\author{Roey Merchav\thanks{Faculty of Data and Decision Science, Technion -Israel Institute of Technology, Haifa 3200003, Israel. Email: merhav.roey@campus.technion.ac.il.} \and Shoham Sabach\stepcounter{footnote}\thanks{School of Operations Research and Information Engineering, Cornell University, Ithaca 14853, USA. Email: ssabach@cornell.edu.}}
\date{}

\begin{document}

\maketitle

\begin{abstract}
    We consider simple bi-level optimization with a nonconvex composite outer objective and a convex composite inner problem. We propose MABPG, a single-loop proximal gradient method that combines dynamic regularization with accelerated inner steps and requires no linear minimization oracle. Using a normalized Frank-Wolfe gap, we establish simultaneous $O(k^{-1/3})$ bounds for the inner objective gap and outer stationarity without compactness or an error-bound assumption. Under a H\"olderian error bound of order $p > 1$ and coercivity of at least one objective, we further obtain $O(k^{-1})$ for the inner objective gap and $O(k^{-1/(2p)})$ for the outer gradient mapping norm. Under quadratic growth, a projection phase variant of MABPG further improves the inner bound to $O(k^{-3})$ while preserving the outer gradient mapping norm bound $O(k^{-1/4})$. Numerical experiments demonstrate the practical effectiveness of MABPG and its advantages over existing methods.
\end{abstract}

\section{Introduction}
    Simple bi-level optimization is a structured subclass of general bi-level optimization \cite{S1952-B,DZ2020} in which an outer objective is minimized over the solution set of an inner optimization problem. This model is useful when the inner problem has multiple minimizers and a secondary criterion selects among them, as in underdetermined inverse problems and over-parameterized regression. In this paper, we consider simple bi-level problems of the form
    \begin{equation*}
        \min_{\bx \in X^{\ast}} \left\{ \omega\left(\bx\right) := \sigma\left(\bx\right) + \psi\left(\bx\right) \right\}, 
    \end{equation*}
    where the constraint set $X^{\ast}$ is given as the optimal solution set of another optimization problem
    \begin{equation*}
        \min_{\bx \in \real^{n}} \left\{ \varphi\left(\bx\right) := f\left(\bx\right) + g\left(\bx\right) \right\}. 
    \end{equation*}
    Here, $f , \sigma : \real^{n} \rightarrow \real$ have Lipschitz continuous gradients, $f$ is convex, $g : \real^{n} \rightarrow \erl$ is proper, lower semicontinuous, and convex, and $\psi : \real^{n} \rightarrow \real$ is convex. The functions $g$ and $\psi$ may include nonsmooth convex regularizers, while $g$ may also encode explicit convex constraints. {\em The smooth outer component $\sigma$ is not assumed to be convex}.

    The outer objective allows us to select an inner solution according to an additional criterion, which need not be convex. Examples include smooth non-convex regularizers, such as the Cauchy regularizer used in sparse representation learning \cite{MKHA2021}. The composite structure of the outer objective function $\omega$ also allows accommodating nonsmooth nonconvex regularizers, such as SCAD and MCP, which can be written as the sum of a convex $\ell_{1}$ term and a smooth concave function with a Lipschitz continuous gradient \cite{WLZ2014}. Another example arises in fair classification, where the squared covariance between predicted probabilities and a sensitive attribute can serve as a smooth nonconvex measure of dependence \cite{JAMH2023}. 

    Viewed as minimizing the non-convex composite function $\omega$ over the closed and convex set $X^{\ast}$, the problem presents a basic difficulty: $X^{\ast}$ is defined implicitly through the minimizers of $\varphi$. Consequently, projections and linear minimization oracles over $X^{\ast}$ need not be available, even when the corresponding operations for the inner problem are easy to compute, and maintaining feasibility would require solving the inner problem exactly. A meaningful approximate stationarity criterion must therefore accommodate points outside $X^{\ast}$, combine outer stationarity with inner objective accuracy, and avoid assuming that $X^{\ast}$ is compact. 

    Simple bi-level optimization has been studied extensively when both $\varphi$ and $\omega$ are convex. Approaches include iterative regularization, also referred to as Tikhonov regularized schemes \cite{TA77-B,FT2008,S07,S2007,GRV2018,AY2019,KY2021,SBY2024,LTVP2025,BCCH2026,MST2026,ST2026}, sequential averaging, viscosity, and bi-level subgradient methods \cite{Xu2004,HS2017,SS2017,SVZ2021,MS2023}, as well as other first-order, cutting-plane, level-set, value-function, penalty, and functional-constraint methods \cite{BS2014,SNK2022,DS2023,JAMH2023,WSJ2024,CJYM2024,ZCXZ2024,CSJW2024,JSSW2026}. These works provide asymptotic convergence results and non-asymptotic bounds under various structural assumptions. However, guarantees for the outer objective gap that rely on convexity do not directly apply when the outer objective is nonconvex.

    For nonconvex outer objectives, the relevant non-asymptotic results can be grouped into three main approaches. The first consists of projection-free conditional gradient methods, including CG-BiO \cite{JAMH2023}, its stochastic extension \cite{CJAHM2023}, and the iteratively regularized method of \cite{GSH2025}. These methods establish bounds involving the inner objective gap and a Frank-Wolfe-type outer stationarity gap, with additional guarantees under a H\"olderian error bound in \cite{JAMH2023}. Their analyses consider smooth objectives and an explicitly given compact convex set over which the inner objective is minimized, referred to here as the inner constraint set. The cutting-plane methods in \cite{JAMH2023,CJAHM2023} use an approximate inner solution for initialization and require linear minimization over the intersection of the inner constraint set with an iteration-dependent halfspace. In contrast, the method in \cite{GSH2025} uses iterative regularization and requires only a standard LMO over the inner constraint set. In this setting, which differs substantially from ours, CG-BiO achieves a faster simultaneous complexity bound than the one we obtain without an error bound.

    A second approach uses gradient directions that balance improvement in the outer objective with progress in the inner problem. Examples include dynamic barrier gradient descent (DBGD) \cite{GLL2021} and gradient-surgery methods \cite{HTNKCA2024}. These methods address smooth problems without explicit constraints and allow the inner objective to be nonconvex. More recently, \cite{CJHM2025} established discrete-time non-asymptotic guarantees for a DBGD-type method using a stationarity criterion defined through the gradients of the two objectives. Such guarantees concern a different notion of stationarity from the one considered here: without further assumptions, they do not provide an inner objective-gap bound or a stationarity estimate for minimizing $\omega$ over the global inner solution set $X^{\ast}$.

    A third approach approximates projections onto the inner solution set. Related inexact projection methods were studied for optimization over variational inequality solution sets in \cite{SY2025}. For a composite convex inner problem and a smooth nonconvex outer objective, \cite{SBY2025} proposed $\text{IPR-VFISTA}_{\text{nc}}$. This method estimates the gradient mapping of the outer objective constrained to $X^{\ast}$ and is particularly relevant to our setting. It has a double-loop structure, with a growing number of accelerated regularized proximal gradient steps used to approximate a projection at each outer iteration. Its non-asymptotic analysis assumes quadratic growth of the inner objective and boundedness of $X^{\ast}$, and its outer objective is smooth. Our setting also allows a nonsmooth convex outer component $\psi$.

    We propose the \emph{Monotone Accelerated Bi-level Proximal Gradient (MABPG)} method, a single-loop algorithm that combines proximal gradient steps for a dynamically regularized objective with FISTA steps for the inner problem \cite{BT2009} through a monotone selection rule in the spirit of \cite{BT2009:1}. Its updates require neither a fixed iteration budget nor a target accuracy. MABPG accommodates nonsmooth composite objectives at both levels and requires no linear minimization oracles or nested loops. Our main contributions are:
    \begin{itemize}[leftmargin=*]
        \item \emph{Simultaneous convergence rates without compactness or an error-bound assumption.} We introduce a normalized Frank-Wolfe gap measure that applies to possibly unbounded inner solution sets and establish its relation to the outer gradient mapping in Section \ref{Sec:StatMeasures}. A suitable choice of the regularization sequence gives an inner objective gap and a normalized outer stationarity gap of order $O(k^{-1/3})$ at the same selected output. These guarantees require neither compactness of $X^{\ast}$ nor a H\"olderian error bound.
        \item \emph{Improved rates under a H\"olderian error bound.} Under a H\"olderian error bound of order $p \geq 1$ and coercivity of at least one objective, we establish stronger inner objective-gap bounds and rates for the outer gradient mapping norm. For $p > 1$, a suitable regularization exponent gives an inner objective gap of order $O\left(k^{-1}\right)$ and an outer gradient mapping norm of order $O\left(k^{-1/\left(2p\right)}\right)$, in both cases at the same selected output. For $p = 1$, the selected output lies in $X^{\ast}$ after finitely many iterations, and from then on its outer gradient mapping norm is of order $O\left(k^{-\left(1 - \gamma\right)/2}\right)$ for any regularization exponent $\gamma \in \left(0 , 1\right)$.
        \item \emph{A faster inner rate via a projection phase.} Under the assumptions of the preceding item with $p = 2$, we complement MABPG with a projection phase based on FBi-PG \cite{MST2026}. The resulting method, MABPG-P, improves the inner objective-gap bound from $O(k^{-1})$ to $O(k^{-3})$ while preserving the outer gradient mapping norm bound $O(k^{-1/4})$ at the same output. This phase uses a number of additional iterations proportional to the first-phase budget, while the updates of MABPG remain independent of the stopping iteration; see Subsection \ref{SSec:MABPGP}.
    \end{itemize}
    The remainder of the paper is organized as follows. Section \ref{Sec:Model} presents the composite bi-level model, standing assumptions, and proximal tools. Section \ref{Sec:StatMeasures} introduces our stationarity measures and establishes their connections to the bi-level gradient mapping and regularized proximal gradient steps. Section \ref{Sec:MABPG} presents MABPG, establishes its convergence rates with and without H\"olderian error bounds, and analyzes MABPG-P. Selected technical proofs are deferred to the appendices. Finally, Section \ref{Sec:Numerics} presents numerical experiments and comparisons with existing methods.

\section{Optimization Model and Mathematical Toolkit} \label{Sec:Model}
    In this paper, we focus on simple bi-level optimization problems with the classical inner-outer structure. Specifically, we consider an \emph{outer optimization problem} of the composite form
    \begin{equation} \label{Prob:OP} \tag{OP}
        \min_{\bx \in X^{\ast}} \left\{ \omega\left(\bx\right) \equiv \sigma\left(\bx\right) + \psi\left(\bx\right) \right\},
    \end{equation}
    where $X^{\ast}$ denotes the set of minimizers of the \emph{inner level problem}, which is itself given by the convex composite model
    \begin{equation} \label{Prob:P} \tag{P}
        \min_{\bx \in \real^{n}} \left\{ \varphi\left(\bx\right) \equiv f\left(\bx\right) + g\left(\bx\right) \right\}.
    \end{equation}
    Throughout the paper, we impose the following standing assumption on the data of the problems under consideration.
    \begin{assumption} \label{A:AssumptionA}
        \begin{itemize}
            \item[$\rm{(i)}$] The functions $f,  \sigma : \real^{n} \rightarrow \real$ are continuously differentiable with Lipschitz continuous gradients, with positive parameters $\beta_{f}$ and $\beta_{\sigma}$, respectively.
            \item[$\rm{(ii)}$] The function $f$ is convex.
            \item[$\rm{(iii)}$] The function $g : \real^{n} \rightarrow \erl$ is proper, lower semicontinuous, and convex. The function $\psi : \real^{n} \rightarrow \real$  is convex.
            \item[$\rm{(iv)}$] The optimal solution set $X^{\ast}$ of problem \eqref{Prob:P} is non-empty.
            \item[$\rm{(v)}$] $\inf_{\bx \in \dom{g}} \omega\left(\bx\right) \equiv \omega^{\ast} > -\infty$.
        \end{itemize}
    \end{assumption}
    For simplicity, we use the notation $\beta \equiv \beta_{f} + \beta_{\sigma}$. For each $k \in \nn$, we define the decay parameter
    \begin{equation} \label{D:Etak}
        \eta_{k} = \left(k + 1\right)^{-\gamma}, \qquad \gamma > 0,
    \end{equation}
    which will govern the dynamic regularization. Accordingly, we define the smooth function $f_{k} \equiv f + \eta_{k}\sigma$, the non-smooth function $g_{k} \equiv g + \eta_{k}\psi$, and the resulting regularized objective
    \begin{equation} \label{D:RegularizedF}
        F_{k} \equiv \varphi + \eta_{k}\omega = f_{k} + g_{k}.
    \end{equation}
    Assumption \ref{A:AssumptionA}(i) guarantees that the gradients $\nabla f_{k}$, $k \in \nn$, are Lipschitz continuous with a uniform Lipschitz constant $\beta$. The following remark makes this observation explicit.
    \begin{remark} \label{R:Lipfk}
        Under Assumption \ref{A:AssumptionA}(i), $f_{k}$ has a Lipschitz continuous gradient with parameter $\beta_{f} + \eta_{k}\beta_{\sigma}$. Since $\eta_{k} \leq 1$ for all $k \in \nn$, we have
        \begin{equation*}
        \beta_{f} + \eta_{k}\beta_{\sigma} \leq \beta_{f} + \beta_{\sigma} = \beta.
        \end{equation*}
        Consequently, $\nabla f_{k}$ is Lipschitz continuous with the uniform constant $\beta$ for all $k \in \nn$.
    \end{remark}
    Requiring $\psi$ to be real-valued may initially appear restrictive. However, when $\psi$ is the indicator function of a constraint set, this requirement holds without loss of generality, as shown in the following remark.
    \begin{remark} \label{R:AssFullDomain}
        Suppose that $\psi \equiv \delta_{C}$ for a nonempty, closed and convex set $C \subseteq \real^{n}$ with $X^{\ast} \cap C \neq \emptyset$. Replacing $g$ by $\tilde{g} \equiv g + \delta_{C}$ and $\psi$ by $\tilde{\psi} \equiv 0$ leaves Problem \eqref{Prob:OP} unchanged. Indeed, since $X^{\ast} \cap C \neq \emptyset$, the set of minimizers of $f + \tilde{g}$ is exactly $X^{\ast} \cap C$, and minimizing $\sigma$ over this set is Problem \eqref{Prob:OP}. Moreover, the new data satisfy Assumption \ref{A:AssumptionA}: $\tilde{\psi} \equiv 0$ is real-valued, and item (v) holds since $\dom{\tilde{g}} \subseteq \dom{g} \cap C$, on which $\sigma = \omega$ is bounded from below.
    \end{remark}    
    Since both inner and outer optimization problems have a composite structure, proximal, proximal gradient, and gradient mappings will play a central role throughout our exposition. We therefore briefly recall the relevant definitions and properties below.

\subsection{Optimization Tools for Composite Models}
    To streamline the presentation, and since the tools introduced in this subsection will be used for both the inner and outer problems, we consider the following general composite optimization model:
     \begin{equation} \label{Composite}
         \min_{\bx \in \real^{n}} \left\{ Q\left(\bx\right) := q_{1}\left(\bx\right) + q_{2}\left(\bx\right) \right\},
     \end{equation}  
     where $q_{1} : \real^{n} \rightarrow \real$ is continuously differentiable, and $q_{2} : \real^{n} \rightarrow \erl$ is proper, lsc, and convex.

    For any $\bx \in \real^{n}$ and a step size $t > 0$, the \emph{proximal gradient mapping} is defined by:
    \begin{equation} \label{D:ProximalGradientOperator}
         T_{t}^{q_{1} , q_{2}}\left(\bx\right) = \argmin_{\bz \in \real^{n}} \left\{ q_{2}\left(\bz\right) + \nabla q_{1}\left(\bx\right)^{T}\bz + \frac{1}{2t}\norm{\bz - \bx}^{2} \right\}.
    \end{equation}
    When $q_{1} \equiv 0$, the proximal gradient mapping $T_{t}^{0 , q_{2}}$ reduces to the (Moreau) proximal mapping of the function $tq_{2}$. When, in addition, $q_{2} \equiv \delta_{C}$ is the indicator function of a nonempty, closed, and convex set $C \subseteq \real^{n}$, it further reduces to the orthogonal projection onto $C$, denoted by $P_{C}\left(\bx\right) := \argmin_{\bz \in C} \norm{\bx - \bz}$.
\medskip

   The associated {\emph gradient mapping} is defined by
    \begin{equation} \label{D:GradientMapping}
        G_{t}^{q_{1} , q_{2}}\left(\bx\right) = t^{-1}\left(\bx - T_{t}^{q_{1} , q_{2}}\left(\bx\right)\right).
    \end{equation}    
    The gradient mapping generalizes the standard gradient to composite optimization models and provides a natural measure of stationarity. In particular, when $q_{2} \equiv 0$, we have $G_{t}^{q_{1} , q_{2}}\left(\bx\right) = \nabla q_{1}\left(\bx\right)$. 
\medskip

    Throughout the remainder of this subsection, we operate under the composite model \eqref{Composite}, where $q_{1} : \real^{n} \rightarrow \real$ is a continuously differentiable function with a Lipschitz continuous gradient with parameter $L$, and $q_{2} : \real^{n} \rightarrow \erl$ is a proper, lsc, and convex function. Recall that $\bx^{\ast}$ is a stationary point of Problem \eqref{Composite} if $-\nabla q_{1}\left(\bx^{\ast}\right) \in \partial q_{2}\left(\bx^{\ast}\right)$. The following standard result characterizes stationarity through the gradient mapping and, in the convex case, characterizes optimality; see \cite[Theorem 10.7 and Corollary 10.8]{B2017-B}.
    \begin{proposition}
        Let $\bx^{\ast} \in \real^{n}$ and $t > 0$. Then
        \begin{itemize}
            \item[$\rm{(i)}$] $G_{t}^{q_{1} , q_{2}}\left(\bx^{\ast}\right) = \bo$ if and only if $\bx^{\ast}$ is a stationary point of Problem \eqref{Composite}. 
            \item[$\rm{(ii)}$] If $q_{1}$ is convex, $G_{t}^{q_{1} , q_{2}}\left(\bx^{\ast}\right) = \bo$ if and only if $\bx^{\ast}$ is an optimal solution to Problem \eqref{Composite}.
        \end{itemize}
    \end{proposition}
    Thus, $\norm{G_{t}^{q_{1},q_{2}}\left(\bx\right)}$ provides a natural stationarity measure for the composite problem, and, when $q_{1}$ is convex, an optimality measure. Moreover, this measure is readily available in algorithms based on the proximal gradient step, since evaluating $G_{t}^{q_{1},q_{2}}$ requires no additional proximal computation beyond that required to obtain $T_{t}^{q_{1},q_{2}}$.

    For later use, we record two additional properties of the proximal gradient mapping. See \cite[Lemma 10.4]{B2017-B}.
    \begin{proposition} \label{P:ProxGradFundProperty}
        Let $\bx \in \real^{n}$ and $0 < t \leq 1/L$. Then
        \begin{itemize}
            \item[$\rm{(i)}$] $G_{t}^{q_{1} , q_{2}}\left(\bx\right) \in \partial q_{2}\left(\bx - tG_{t}^{q_{1}, q_{2}}\left(\bx\right)\right) + \nabla q_{1}\left(\bx\right)$. 
            \item[$\rm{(ii)}$] $Q\left(T_{t}^{q_{1} , q_{2}}\left(\bx\right)\right) \leq Q\left(\bx\right) - \frac{1}{2t}\norm{\bx - T_{t}^{q_{1} , q_{2}}\left(\bx\right)}^{2}$.
        \end{itemize}
    \end{proposition}    
    Finally, since the proposed method employs accelerated schemes for composite optimization problems, we will also use several properties of FISTA \cite{BT2009}, which we record next.
    \begin{proposition}[FISTA properties] \label{P:FISTAprop}
        Let $\Sequp{\bz}{k}$ be a sequence generated by FISTA applied on Problem \eqref{Composite} with initialization $\bz^{0} \in \real^{n}$ and a fixed step size $t = 1/L$. Assume that the function $q_{1}$ is convex. Then, for all $\bx^{\ast} \in \argmin_{\bx \in \real^{n}} Q\left(\bx\right)$ and $k \geq 1$, the following items hold true:
        \begin{itemize}
            \item[$\rm{(i)}$] $Q\left(\bz^{k}\right) - Q\left(\bx^{\ast}\right) \leq 2L\norm{\bz^{0} - \bx^{\ast}}^{2}\left(k + 1\right)^{-2}$.
            \item[$\rm{(ii)}$] $\norm{\bz^{k} - \bx^{\ast}} \leq \norm{\bz^{0} - \bx^{\ast}}$. In particular, the sequence $\Sequp{\bz}{k}$ is bounded.
        \end{itemize}
    \end{proposition}
    \begin{proof}
        Item (i) is proved in \cite[Theorem 4.4]{BT2009}. For item (ii), from the proof of \cite{BT2009}, we also have that
        \begin{equation*}
            \frac{2}{L}t_{k - 1}^{2}\left(Q\left(\bz^{k}\right) - Q\left(\bx^{\ast}\right)\right) + \norm{\bw^{k} - \bx^{\ast}}^{2} \leq \norm{\bz^{0} - \bx^{\ast}}^{2},
        \end{equation*}
        where $\Sequp{\bw}{k}$ is the FISTA auxiliary sequence defined by $\bw^{k} = \bz^{k - 1} + t_{k - 1}\left(\bz^{k} - \bz^{k - 1}\right)$. Since $\bx^{\ast}$ is a minimizer of $Q$, we have that $Q\left(\bz^{k}\right) \geq Q\left(\bx^{\ast}\right)$, which implies
        \begin{equation} \label{P:FISTAprop:1}
            \norm{\bw^{k} - \bx^{\ast}} \leq \norm{\bz^{0} - \bx^{\ast}}.
        \end{equation}
        Rearranging $\bw^{k}$ gives $\bz^{k} = \left(1 - t_{k - 1}^{-1}\right)\bz^{k - 1} + t_{k - 1}^{-1}\bw^{k}$. Since $t_{k - 1} \geq 1$ it follows that $t_{k - 1}^{-1} \in \left(0 , 1\right]$, and thus $\bz^{k}$ is a convex combination of $\bz^{k - 1}$ and $\bw^{k}$. By the convexity of the norm, we obtain from \eqref{P:FISTAprop:1} that
        \begin{equation*}
            \norm{\bz^{k} - \bx^{\ast}} \leq \left(1 - t_{k - 1}^{-1}\right)\norm{\bz^{k - 1} - \bx^{\ast}} + t_{k - 1}^{-1}\norm{\bw^{k} - \bx^{\ast}} \leq \left(1 - t_{k - 1}^{-1}\right)\norm{\bz^{k - 1} - \bx^{\ast}} + t_{k - 1}^{-1}\norm{\bz^{0} - \bx^{\ast}}.
        \end{equation*}
        The result now follows by a simple induction argument.
    \end{proof}

\section{Stationarity in Simple Bi-level Optimization} \label{Sec:StatMeasures}
    In this paper, we consider simple bi-level optimization problems in which the outer objective function is not necessarily convex; recall that no convexity assumption is imposed on $\sigma$. This nonconvexity necessitates a careful choice of an appropriate stationarity measure. To this end, we view the bi-level problem \eqref{Prob:OP} as the constrained optimization problem
   \begin{equation*}
        \min_{\bx \in X^{\ast}} \left\{ \omega\left(\bx\right) \equiv \sigma\left(\bx\right) + \psi\left(\bx\right) \right\},
    \end{equation*}
    where $X^{\ast}$ denotes the \emph{convex} set of minimizers of the inner level problem \eqref{Prob:P}. It follows that $\bx$ is a stationary point of Problem \eqref{Prob:OP} if $\bo \in \nabla \sigma\left(\bx\right) + \partial \psi\left(\bx\right) + \mathcal{N}_{X^{\ast}}\left(\bx\right)$, where $\mathcal{N}_{X^{\ast}}\left(\bx\right)$ is the normal cone to $X^{\ast}$ at $\bx$

    Problem \eqref{Prob:OP} contains two sources of nonsmoothness: the function $\psi$ and the constraint set $X^{\ast}$. Accordingly, there are two natural ways to characterize stationarity. The first is through the gradient mapping introduced in \eqref{D:GradientMapping}, applied to the composite function $q_{1} \equiv \sigma$ and $q_{2} \equiv \psi + \delta_{X^{\ast}}$. The second is through the generalized Frank-Wolfe (FW) gap introduced in \cite{B2015}; see also \cite{ZZ2023}. The latter avoids incorporating the indicator function $\delta_{X^{\ast}}$ into the nonsmooth component and is defined as follows.
    \begin{definition}[Generalized Frank-Wolfe gap] \label{D:GenFWGap}
        The generalized Frank-Wolfe gap for Problem \eqref{Prob:OP} at any point $\bx \in X^{\ast}$ is defined as:
        \begin{equation}
            {\hat \GGG}_{OP}\left(\bx\right) \equiv \sup_{\by \in X^{\ast}} \left\{ \act{\nabla \sigma\left(\bx\right) , \bx - \by} + \psi\left(\bx\right) - \psi\left(\by\right) \right\}.
        \end{equation}
    \end{definition}
    The generalized FW gap is a useful measure for stationarity of Problem \eqref{Prob:OP} due to the following basic properties.
    \begin{proposition}
        \begin{itemize}
            \item[$\rm{(i)}$] ${\hat \GGG}_{OP}\left(\bx\right) \geq 0$ for all $\bx \in X^{\ast}$.
            \item[$\rm{(ii)}$] ${\hat \GGG}_{OP}\left(\bx^{\ast}\right) = 0$ if and only if $\bx^{\ast}$ is a stationary point of Problem \eqref{Prob:OP}. 
            \item[$\rm{(iii)}$] If $\sigma$ is convex, ${\hat \GGG}_{OP}\left(\bx^{\ast}\right) = 0$ if and only if $\bx^{\ast}$ is an optimal solution to Problem \eqref{Prob:OP}.
        \end{itemize}
    \end{proposition}
    The generalized FW gap is widely used in Linear Oracle (LO)-based optimization methods \cite{B2015,ZZ2023} and has also been employed in the context of simple bi-level optimization. For example, \cite{JAMH2023} uses this measure when the outer objective function is smooth, corresponding to $\psi \equiv 0$ in our notation. In the LO setting, however, the generalized FW gap is typically considered over a compact constraint set. In the bi-level setting, this amounts to requiring the inner-level solution set $X^{\ast}$ to be compact. In \cite{JAMH2023}, compactness is ensured by assuming that the inner-level nonsmooth function is the indicator function of a compact set. In this paper, we develop a modification of the generalized FW gap that remains suitable for simple bi-level optimization problems even when the inner-level solution set $X^{\ast}$ is not compact.
    \begin{definition}[Normalized Frank-Wolfe gap] \label{D:NormFWGap}
        The normalized Frank-Wolfe gap for Problem \eqref{Prob:OP} at any $\bx \in \real^{n}$ is defined as:
        \begin{equation}
            \GGG_{OP}\left(\bx\right) \equiv \sup_{\by \in X^{\ast} \setminus \{ \bx \}} \frac{\act{\nabla \sigma\left(\bx\right) , \bx - \by} + \psi\left(\bx\right) - \psi\left(\by\right)}{\norm{\bx - \by}},
        \end{equation}
        where we use the convention that $\sup \emptyset = -\infty$.
    \end{definition}
    The FW gap measures are well-suited to standard LO-based methods, in which the iterates remain feasible with respect to the constraint set. However, algorithms for simple bi-level optimization generally generate iterates that do not belong to $X^{\ast}$; this is also the case for the method proposed in this paper. This means the normalized FW gap alone does not capture progress toward feasibility. Motivated by \cite{JAMH2023}, we therefore combine it with the optimality gap of the inner problem, leading to the following bi-level stationarity measure.
    \begin{definition}[$(\epsilon_{\varphi} , \epsilon_{\omega})$-stationary point] \label{D:EpsStat}
        A point $\bx \in \real^{n}$ is called an $(\epsilon_{\varphi} , \epsilon_{\omega})$-stationary point of Problem \eqref{Prob:OP} if it satisfies:
        \begin{enumerate}
            \item[$\rm{(i)}$] $\varphi\left(\bx\right) - \min_{\by \in \real^{n}} \varphi\left(\by\right) \leq \epsilon_{\varphi}$,
            \item[$\rm{(ii)}$] $\GGG_{OP}\left(\bx\right) \leq \epsilon_{\omega}$.
        \end{enumerate}
    \end{definition}
    To further justify this measure, we next show that, for convex simple bi-level optimization problems, it characterizes optimal solutions of Problem \eqref{Prob:OP}.
    \begin{proposition} \label{P:CFWValidity}
        Suppose that Assumption \ref{A:AssumptionA} holds, and the function $\omega$ is convex. Then, a point $\bx \in \real^{n}$ is an optimal solution to Problem \eqref{Prob:OP} if and only if it is a $(0,  0)$-stationary point. 
    \end{proposition} 
    \begin{proof}
        A point $\bx$ is $\left(0 , 0\right)$-stationary if and only if $\bx \in X^{\ast}$ and
        \begin{equation*}
            \act{\nabla \sigma\left(\bx\right) , \by - \bx} + \psi\left(\by\right) - \psi\left(\bx\right) \geq 0 \quad \text{for all } \by \in X^{\ast}.
        \end{equation*}
        This is equivalent to $\bx$ minimizing the convex function $\act{\nabla \sigma\left(\bx\right) , \cdot} + \psi + \delta_{X^{\ast}}$. Since $\dom{\psi} = \real^{n}$, the optimality condition and the subdifferential sum rule give
        \begin{equation*}
            \bo \in \nabla \sigma\left(\bx\right) + \partial \psi\left(\bx\right) + \mathcal{N}_{X^{\ast}}\left(\bx\right) = \partial \omega\left(\bx\right) + \mathcal{N}_{X^{\ast}}\left(\bx\right).
        \end{equation*}
        By the convexity of $\omega$, this condition is equivalent to $\bx$ being an optimal solution to Problem \eqref{Prob:OP}.
    \end{proof}
    The next result bounds the outer gradient mapping by the positive part of the normalized FW gap and the distance to $X^{\ast}$. Its proof is given in Appendix \ref{App:StationarityBounds}.
    \begin{proposition} \label{P:GradMapFWgap}
        For any $\bx \in \real^{n}$ and each step size $t > 0$, we have
        \begin{equation} \label{P:GradMapFWgap:1}
            \norm{G_{t}^{\sigma , \psi + \delta_{X^{\ast}}}\left(\bx\right)} \leq 2\pos{\GGG_{OP}\left(\bx\right)} + \frac{1}{t}\norm{\bx^{\ast} - \bx} + \sqrt{\frac{2}{t}\norm{\nabla \sigma\left(\bx\right) + \xi} \cdot \norm{\bx^{\ast} - \bx}},
        \end{equation}
        where $\bx^{\ast} = P_{X^{\ast}}\left(\bx\right)$, $\xi \in \partial \psi\left(\bx^{\ast}\right)$, and $\pos{s} := \max\left\{s , 0\right\}$ for $s \in \left[-\infty , +\infty\right]$.
    \end{proposition}
    The next result shows that a proximal gradient step on the regularized objective controls the normalized FW gap. Its proof is also given in Appendix \ref{App:StationarityBounds}.
    \begin{proposition} \label{P:RegStepFW}
        Let $\eta > 0$, $f_{\eta} := f + \eta\sigma$ and $g_{\eta} := g + \eta\psi$, and let $\beta > 0$ be a Lipschitz constant of $\nabla f_{\eta}$. For $\bv \in \real^{n}$, we denote $\bx^{+} := T_{1/\beta}^{f_{\eta} , g_{\eta}}\left(\bv\right)$ and $G_{1/\beta}\left(\bv\right) := G_{1/\beta}^{f_{\eta} , g_{\eta}}\left(\bv\right)$. Then, for any $\bx \in \real^{n}$, we have
        \begin{equation} \label{P:RegStepFW:1}
            \varphi\left(\bx^{+}\right) - \varphi\left(\bx\right) + \eta\left(\psi\left(\bx^{+}\right) - \psi\left(\bx\right) + \act{\nabla \sigma\left(\bx^{+}\right) , \bx^{+} - \bx}\right) \leq 2\norm{G_{1/\beta}\left(\bv\right)} \cdot \norm{\bx^{+} - \bx}.
        \end{equation}
        In particular, $\GGG_{OP}\left(\bx^{+}\right) \leq 2\eta^{-1}\norm{G_{1/\beta}\left(\bv\right)}$.
    \end{proposition}
    
\section{Monotone Accelerated Bi-level Proximal Gradient Algorithm} \label{Sec:MABPG}
    We now introduce the Monotone Accelerated Bi-level Proximal Gradient (MABPG) algorithm for Problem \eqref{Prob:OP}. Recall the regularized objective $F_{k} = \varphi + \eta_{k}\omega$ defined in \eqref{D:RegularizedF}, where $\eta_{k} = \left(k + 1\right)^{-\gamma}$ and $\gamma > 0$. 

    \noindent\fbox{\parbox{\textwidth}{
    {\bf Monotone Accelerated Bi-level Proximal Gradient (MABPG)}
    \begin{itemize}[leftmargin=*]
        \item[] {\bf Input:} Lipschitz constant $\beta > 0$. Regularization parameter $\gamma > 0$.
        \item[] {\bf Initialization}: $\bz^{-1} = \bz^{0} = \bv^{-1} = \bv^{0} \in \real^{n}$, $t_{-1} = 1$, $t_{0} = 1$. 
        \item[] {\bf General Step} ($k = 0, 1, 2 , \ldots$):
            \begin{itemize}
                \item[] Let $\eta_{k} = \left(k + 1\right)^{-\gamma}$. A Proximal Gradient step on the regularized problem
                    \begin{equation} \label{Bi-AGUpdate:3}
                        \bx^{k + 1} = T_{1/\beta}^{f_{k} , g_{k}}\left(\bv^{k}\right).
                    \end{equation}
                \item[] A FISTA step on the inner problem
                    \begin{equation} \label{Bi-AGUpdate:1}
                        \bz^{k + 1} = T_{1/\beta}^{f , g}\left(\bz^{k} + \tfrac{t_{k - 1} - 1}{t_{k}}\left(\bz^{k} - \bz^{k - 1}\right)\right).
                    \end{equation}
                \item[] Main step
                \begin{equation} \label{Bi-AGUpdate:5}
                    \bv^{k + 1} \in \argmin \left\{ F_{k}\left(\bx\right) : \, \bx \in \{ \bz^{k + 1} , \bx^{k + 1} \} \right\}.
                \end{equation}
                Updating the parameters $t_{k + 1} = \tfrac{1 + \sqrt{1 + 4t_{k}^{2}}}{2}$.
            \end{itemize}
    \end{itemize}}}
\vspace{0.15in}

    At each iteration, MABPG generates two candidate points. The point $\bz^{k + 1}$ is obtained by applying FISTA to the inner problem and provides an accelerated inner objective bound. The point $\bx^{k + 1}$ is obtained by a proximal gradient step on the regularized problem and provides the descent estimate used in the analysis. The next iterate is selected among these candidates by minimizing $F_{k}$.

    \begin{remark}[Additional candidates] \label{R:AddCand}
        The convergence analysis only requires that the selected iterate satisfy
        \begin{equation*}
            F_{k}\left(\bv^{k + 1}\right) \leq \min \left\{ F_{k}\left(\bx^{k + 1}\right) , F_{k}\left(\bz^{k + 1}\right) \right\}.
        \end{equation*}
        Consequently, additional candidate points may be included in the comparison step without changing the convergence guarantees. In the numerical experiments, we use the additional candidate
        \begin{equation*}
            \bu^{k + 1} = T_{1/\beta}^{f_{k} , g_{k}}\left(\bv^{k} + \frac{t_{k - 1} - 1}{t_{k}}\left(\bv^{k} - \bv^{k - 1}\right)\right),
        \end{equation*}
        and select $\bv^{k + 1}$ by minimizing $F_{k}$ over $\left\{ \bx^{k + 1} , \bz^{k + 1} , \bu^{k + 1} \right\}$.
    \end{remark}
    \begin{remark}[Well-definedness] \label{R:Welld}
        Let $k \geq 1$. All vectors $\bz^{k}$, $\bu^{k}$ and $\bx^{k}$ are outputs of proximal gradient steps with respect to $g$ or $g_{k - 1} = g + \eta_{k -1}\psi$. Since $\psi$ is a real-valued function (see Assumption \ref{A:AssumptionA}(iii)), we have that $\dom{g_{k}} = \dom{g}$, and thus $\bv^{k} \in \dom{g}$. 
    \end{remark}
    Our theoretical analysis is organized into two complementary stages, which together provide a comprehensive convergence analysis of the proposed algorithm. In the first stage, we focus on the inner-level problem and establish bounds for the dynamically regularized objectives $F_{k}$, which we then use to derive a convergence rate for the inner objective $\varphi$. In the second stage, we analyze the convergence of the normalized FW gap $\GGG_{OP}$. We first establish a convergence rate for the norm of the gradient mapping $G_{1/\beta}^{f_{k} , g_{k}}$ associated with the regularized problem, and then combine this result with the bounds obtained in the first stage to derive a convergence rate in terms of $(\epsilon_{\varphi} , \epsilon_{\omega})$-stationarity; see Definition \ref{D:EpsStat}.

\subsection{Inner Convergence Rate of MABPG Algorithm}
    The analysis of MABPG is driven by a single quantity, which combines the inner objective gap with the regularized outer objective. Throughout, we denote by $\varphi^{\ast} := \min_{\bx \in \real^{n}} \varphi\left(\bx\right)$ the optimal value of the inner problem and by $d_{0} := \dist\left(\bv^{0} , X^{\ast}\right)$ the distance of the starting point to $X^{\ast}$. Moreover, we define
    \begin{equation*}
        \Omega := \sup_{k \in \nn} \omega\left(\bz^{k}\right).
    \end{equation*}
    We claim that $\Omega < \infty$. Indeed, Proposition \ref{P:FISTAprop}(ii) implies that $\Sequp{\bz}{k}$ is bounded, while Assumption \ref{A:AssumptionA}(iii) then guarantees that $\left\{ \omega\left(\bz^{k}\right) \right\}_{k \in \nn}$ is bounded. Hence, $\Omega < \infty$.

    Let $\left\{ \left(\bx^{k} , \bv^{k} \right) \right\}_{k \in \nn}$ be a sequence generated by MABPG. For all $k \geq 1$, we define
    \begin{equation*}
        E_{k} := \varphi\left(\bv^{k}\right) - \varphi^{\ast} + \eta_{k - 1}\left(\omega\left(\bv^{k}\right) - \omega^{\ast}\right).
    \end{equation*}
    Notice that $E_{k} = F_{k - 1}\left(\bv^{k}\right) - \varphi^{\ast} - \eta_{k - 1}\omega^{\ast}$. The next lemma collects its three basic properties.
    \begin{lemma} \label{L:Energy}
        Let $\left\{ \left(\bx^{k} , \bv^{k} \right) \right\}_{k \in \nn}$ be a sequence generated by MABPG with $\gamma > 0$. Then, for all $k \geq 1$, the following items hold true.
        \begin{itemize}
            \item[$\rm{(i)}$] $E_{k} \geq 0$.
            \item[$\rm{(ii)}$] $E_{k} \leq 2\beta d_{0}^{2}\left(k + 1\right)^{-2} + \eta_{k - 1}\left(\Omega - \omega^{\ast}\right)$.
            \item[$\rm{(iii)}$] $E_{k + 1} \leq \varphi\left(\bx^{k + 1}\right) - \varphi^{\ast} + \eta_{k}\left(\omega\left(\bx^{k + 1}\right) - \omega^{\ast}\right) \leq E_{k} - \dfrac{1}{2\beta}\norm{G_{1/\beta}^{f_{k} , g_{k}}\left(\bv^{k}\right)}^{2}$.
        \end{itemize}
    \end{lemma}
    \begin{proof}
        From Remark \ref{R:Welld}, we have that $\bv^{k} \in \dom{g}$. Therefore, $\varphi\left(\bv^{k}\right) \geq \varphi^{\ast}$ and $\omega\left(\bv^{k}\right) \geq \omega^{\ast}$ (see Assumption \ref{A:AssumptionA}(v)), which proves the first item.

        Let $\bx^{\ast} := P_{X^{\ast}}\left(\bv^{0}\right)$. From Step \eqref{Bi-AGUpdate:5}, we have $F_{k - 1}\left(\bv^{k}\right) \leq F_{k - 1}\left(\bz^{k}\right)$. Hence
        \begin{align*}
            E_{k} = F_{k - 1}\left(\bv^{k}\right) - \varphi^{\ast} - \eta_{k - 1}\omega^{\ast} & \leq F_{k - 1}\left(\bz^{k}\right) - \varphi^{\ast} - \eta_{k - 1}\omega^{\ast} \\
            & = \varphi\left(\bz^{k}\right) - \varphi\left(\bx^{\ast}\right) + \eta_{k - 1}\left(\omega\left(\bz^{k}\right) - \omega^{\ast}\right) \\
            & \leq \frac{2\beta\norm{\bz^{0} - \bx^{\ast}}^{2}}{\left(k + 1\right)^{2}} + \eta_{k - 1}\left(\Omega - \omega^{\ast}\right),
        \end{align*}
        where the last inequality follows from Proposition \ref{P:FISTAprop}(i) and the definition of $\Omega$. Since $\bz^{0} = \bv^{0}$, we have $\norm{\bz^{0} - \bx^{\ast}} = d_{0}$, which proves the second item.

        From Step \eqref{Bi-AGUpdate:3}, we have that $\bx^{k + 1} = T_{1/\beta}^{f_{k} , g_{k}}\left(\bv^{k}\right)$. Since $f_{k}$ is $\beta$-smooth (see Remark \ref{R:Lipfk}), we can apply Proposition \ref{P:ProxGradFundProperty}(ii) with $q_{1} = f_{k}$, $q_{2} = g_{k}$ and $t = 1/\beta$, to obtain that
        \begin{equation} \label{L:Energy:1}
            F_{k}\left(\bx^{k + 1}\right) \leq F_{k}\left(\bv^{k}\right) - \frac{\beta}{2}\norm{\bv^{k} - \bx^{k + 1}}^{2} = F_{k}\left(\bv^{k}\right) - \frac{1}{2\beta}\norm{G_{1/\beta}^{f_{k} , g_{k}}\left(\bv^{k}\right)}^{2},
        \end{equation}
        where the equality follows from the fact that $\bx^{k + 1} = \bv^{k} - \beta^{-1}G_{1/\beta}^{f_{k} , g_{k}}\left(\bv^{k}\right)$. Moreover, from Step \eqref{Bi-AGUpdate:5}, we have $F_{k}\left(\bv^{k + 1}\right) \leq F_{k}\left(\bx^{k + 1}\right)$. Hence
        \begin{equation*}
            E_{k + 1} = F_{k}\left(\bv^{k + 1}\right) - \varphi^{\ast} - \eta_{k}\omega^{\ast} \leq F_{k}\left(\bx^{k + 1}\right) - \varphi^{\ast} - \eta_{k}\omega^{\ast} = \varphi\left(\bx^{k + 1}\right) - \varphi^{\ast} + \eta_{k}\left(\omega\left(\bx^{k + 1}\right) - \omega^{\ast}\right),
        \end{equation*}
        which proves the first inequality of Item (iii). Moreover, using \eqref{L:Energy:1} and after subtracting $\varphi^{\ast} + \eta_{k}\omega^{\ast}$ from both sides, we obtain
        \begin{equation*}
            \varphi\left(\bx^{k + 1}\right) - \varphi^{\ast} + \eta_{k}\left(\omega\left(\bx^{k + 1}\right) - \omega^{\ast}\right) \leq \varphi\left(\bv^{k}\right) - \varphi^{\ast} + \eta_{k}\left(\omega\left(\bv^{k}\right) - \omega^{\ast}\right) - \frac{1}{2\beta}\norm{G_{1/\beta}^{f_{k} , g_{k}}\left(\bv^{k}\right)}^{2}.
        \end{equation*}
        Since $\eta_{k} \leq \eta_{k - 1}$ and $\omega\left(\bv^{k}\right) \geq \omega^{\ast}$, we have $\eta_{k}\left(\omega\left(\bv^{k}\right) - \omega^{\ast}\right) \leq \eta_{k - 1}\left(\omega\left(\bv^{k}\right) - \omega^{\ast}\right)$, and thus the right-hand side is bounded from above by $E_{k} - \left(2\beta\right)^{-1}\norm{G_{1/\beta}^{f_{k} , g_{k}}\left(\bv^{k}\right)}^{2}$. This proves the third item.
    \end{proof}
    The inner rate of MABPG now follows immediately.
    \begin{theorem}[Inner rate] \label{T:InnerRate}
        Let $\left\{ \left(\bx^{k} , \bv^{k}\right) \right\}_{k \in \nn}$ be a sequence generated by MABPG with $\gamma > 0$. Then, for all $k \geq 1$, we have
        \begin{equation*}
            \max\left\{ \varphi\left(\bx^{k + 1}\right) , \varphi\left(\bv^{k + 1}\right) \right\} - \varphi^{\ast} \leq \frac{2\beta d_{0}^{2}}{\left(k + 1\right)^{2}} + \frac{\Omega - \omega^{\ast}}{k^{\gamma}}.
        \end{equation*}
    \end{theorem}
    \begin{proof}
        Since $\omega\left(\bx^{k + 1}\right) \geq \omega^{\ast}$ and $\omega\left(\bv^{k + 1}\right) \geq \omega^{\ast}$, which follows from Remark \ref{R:Welld}, Lemma \ref{L:Energy}(iii) yields $\varphi\left(\bx^{k + 1}\right) - \varphi^{\ast} \leq E_{k}$ and $\varphi\left(\bv^{k + 1}\right) - \varphi^{\ast} \leq E_{k + 1} \leq E_{k}$. The result now follows from Lemma \ref{L:Energy}(ii), since $\eta_{k - 1} = k^{-\gamma}$.
    \end{proof}
    Another immediate consequence is a bound on the outer objective along the iterates.
    \begin{lemma} \label{L:FunctionBound}
        Let $\left\{ \left(\bx^{k} , \bv^{k} \right) \right\}_{k \in \nn}$ be a sequence generated by MABPG with $0 < \gamma \leq 2$. Then, for all $k \geq 1$, we have
        \begin{equation*}
            \omega\left(\bx^{k + 1}\right) - \omega^{\ast} \leq 2\beta d_{0}^{2} + 2^{\gamma}\left(\Omega - \omega^{\ast}\right).
        \end{equation*}
        In particular, the sequences $\left\{ \varphi\left(\bx^{k}\right) \right\}_{k \in \nn}$ and $\left\{ \omega\left(\bx^{k}\right) \right\}_{k \in \nn}$ are bounded.
    \end{lemma}
    \begin{proof}
        Since $\varphi \geq \varphi^{\ast}$, Lemma \ref{L:Energy}(iii) yields $\eta_{k}\left(\omega\left(\bx^{k + 1}\right) - \omega^{\ast}\right) \leq E_{k}$. Dividing by $\eta_{k}$ and using Lemma \ref{L:Energy}(ii), we obtain
        \begin{equation*}
            \omega\left(\bx^{k + 1}\right) - \omega^{\ast} \leq \frac{2\beta d_{0}^{2}}{\left(k + 1\right)^{2 - \gamma}} + \frac{\eta_{k - 1}}{\eta_{k}}\left(\Omega - \omega^{\ast}\right) \leq 2\beta d_{0}^{2} + 2^{\gamma}\left(\Omega - \omega^{\ast}\right),
        \end{equation*}
        where the last inequality follows from the facts that $\left(k + 1\right)^{2 - \gamma} \geq 1$ (recall that $0 < \gamma \leq 2$) and $\eta_{k - 1}\eta_{k}^{-1} = \left(1 + 1/k\right)^{\gamma} \leq 2^{\gamma}$. The boundedness of $\left\{ \varphi\left(\bx^{k}\right) \right\}_{k \in \nn}$ follows from Theorem \ref{T:InnerRate}.
    \end{proof}
    
\subsection{Outer Convergence Rate of MABPG Algorithm}
    The main computational steps of MABPG involve the proximal gradient mapping. Using its relation to the gradient mapping (see \eqref{D:GradientMapping}), these steps can equivalently be expressed in terms of the gradient mapping. For example, Step \eqref{Bi-AGUpdate:3} can be written equivalently as
    \begin{equation} \label{ProxGradStep}
        \bx^{k + 1} = \bv^{k} - \beta^{-1}G_{1/\beta}^{f_{k} , g_{k}}\left(\bv^{k}\right).
    \end{equation}
    An immediate consequence of Lemma \ref{L:Energy} gives us rates on the gradient mapping of the regularized problem. For simplicity, we will use from now on the following notations (recall $d_{0} := \dist\left(\bv^{0} , X^{\ast}\right)$ an $\Omega := \sup_{k \in \nn} \omega\left(\bz^{k}\right)$):
    \begin{equation} \label{Notations}
         C := \sqrt{4\beta^{2}d_{0}^{2} + 2\beta\left(\Omega - \omega^{\ast}\right)}, \quad \text{and} \quad k_{best} \in \argmin_{k \leq s \leq 2k - 1} \norm{G_{1/\beta}^{f_{s} , g_{s}}\left(\bv^{s}\right)}.
    \end{equation}
    \begin{proposition} \label{P:GkBound}
        Let $\left\{\left(\bx^{k} , \bv^{k} \right)\right\}_{k \in \nn}$ be a sequence generated by MABPG with $0 < \gamma \leq 2$. Then, for all $k \geq 1$, we have
        \begin{equation} \label{P:GkBound:1}
            \norm{G_{1/\beta}^{f_{k} , g_{k}}\left(\bv^{k}\right)} \leq \frac{C}{k^{\gamma/2}},
        \end{equation}
        and $\left\{ \norm{G_{1/\beta}^{f_{k} , g_{k}}\left(\bv^{k}\right)} \right\}_{k \in \nn}$ converges to zero. Moreover
        \begin{equation} \label{P:GkBound:2}
            \norm{G_{1/\beta}^{f_{k_{best}} , g_{k_{best}}}\left(\bv^{k_{best}}\right)} \leq \frac{C}{k^{\left(1 + \gamma\right)/2}}.
        \end{equation}
    \end{proposition}
    \begin{proof}
        From Lemma \ref{L:Energy}(i) and (iii), we obtain that
        \begin{equation*}
            \frac{1}{2\beta}\norm{G_{1/\beta}^{f_{k} , g_{k}}\left(\bv^{k}\right)}^{2} \leq E_{k} - E_{k + 1} \leq E_{k} \leq \frac{2\beta d_{0}^{2}}{\left(k + 1\right)^{2}} + \frac{\Omega - \omega^{\ast}}{k^{\gamma}} \leq k^{-\gamma}\left(2\beta d_{0}^{2} + \Omega - \omega^{\ast}\right),
        \end{equation*}
        where the last inequality follows from the facts that $\left(k + 1\right)^{2} \geq \left(k + 1\right)^{2 - \gamma}k^{\gamma}$ and $\left(k + 1\right)^{2 - \gamma} \geq 1$ (recall that $0 < \gamma \leq 2$). Multiplying both sides by $2\beta$ yields the first result. Moreover, summing the first inequality from $s = k$ to $2k - 1$, yields
        \begin{equation*}
            \frac{1}{2\beta}\sum_{s = k}^{2k - 1} \norm{G_{1/\beta}^{f_{s} , g_{s}}\left(\bv^{s}\right)}^{2} \leq E_{k} - E_{2k} \leq E_{k} \leq k^{-\gamma}\left(2\beta d_{0}^{2} + \Omega - \omega^{\ast}\right).
        \end{equation*}
        Using the definition of $k_{best}$, we have
        \begin{equation*}
            k \cdot \norm{G_{1/\beta}^{f_{k_{best}} , g_{k_{best}}}\left(\bv^{k_{best}}\right)}^{2} \leq \sum_{s = k}^{2k - 1} \norm{G_{1/\beta}^{f_{s} , g_{s}}\left(\bv^{s}\right)}^{2},
        \end{equation*}
        and the second result follows using the same arguments.
    \end{proof}
    The selected index $k_{best}$ is therefore used to obtain the sharper residual rate of Proposition \ref{P:GkBound}. For $0<\gamma\leq2$, selection improves the exponent from $\gamma/2$ to $\left(1 + \gamma\right)/2$. This improvement is essential for the outer estimate. Following the discussion in Section \ref{Sec:StatMeasures}, we next present our outer rate of convergence result in terms of the normalized FW Gap.
    \begin{theorem}[Outer rate] \label{T:GkBound}
        Let $\left\{\left(\bx^{k} , \bv^{k} \right)\right\}_{k \in \nn}$ be a sequence generated by MABPG with $0 < \gamma \leq 1$. Then, for any $k \geq 1$, we have
        \begin{equation} \label{T:GkBound:1}
            \GGG_{OP}\left(\bx^{k_{best} + 1}\right) \leq \frac{4C}{k^{\left(1 - \gamma\right)/2}},
        \end{equation}
        with $C$ and $k_{best}$ are as defined in \eqref{Notations}.
    \end{theorem}
    \begin{proof}
        From Propositions \ref{P:GkBound} and \ref{P:RegStepFW} (with $\bx^{+} = \bx^{k_{best} + 1}$, $\bv = \bv^{k_{best}}$ and $\eta = \eta_{k_{best}}$) we obtain that
        \begin{equation*}
            \frac{\eta_{k_{best}}}{2} \cdot \GGG_{OP}\left(\bx^{k_{best} + 1}\right) \leq \norm{G_{1/\beta}^{f_{k_{best}} , g_{k_{best}}}\left(\bv^{k_{best}}\right)} \leq \frac{C}{k^{\left(1 + \gamma\right)/2}}.
        \end{equation*}
        Since $k_{best} \leq 2k - 1$ it follows that $\eta_{k_{best}} \geq \eta_{2k - 1} = \left(2k\right)^{-\gamma} \geq 2^{-1}k^{-\gamma}$ (since $\gamma \leq 1$), we obtain that
        \begin{equation*}
            \GGG_{OP}\left(\bx^{k_{best} + 1}\right) \leq 4k^{\gamma} \cdot \frac{C}{k^{\left(1 + \gamma\right)/2}},
        \end{equation*}
        which proves the desired result.
    \end{proof}
    
\subsection{Simultaneous Convergence}   
    \begin{theorem}[Simultaneous rates] \label{T:SimulConvergenceRate}
        Let $\left\{ \left(\bx^{k} , \bv^{k} \right) \right\}_{k \in \nn}$ be a sequence generated by MABPG with $\gamma = 1/3$. Then, for any $k \geq 1$, the iterate $\bx^{k_{best} + 1}$ is an $(\epsilon_{\varphi} , \epsilon_{\omega})$-stationary point of Problem \eqref{Prob:OP} with
        \begin{equation}
            \epsilon_{\varphi} = \frac{2\beta d_{0}^{2}}{\left(k + 1\right)^{2}} + \frac{\Omega - \omega^{\ast}}{k^{1/3}}, \quad \text{and} \quad \epsilon_{\omega} = \frac{4C}{k^{1/3}}.
        \end{equation} 
        Thus, both the inner sub-optimality and the positive part of the normalized FW gap converge simultaneously at a rate of $\mathcal{O}(k^{-1/3})$.
    \end{theorem}
    \begin{proof}
        The bound for $\epsilon_{\varphi}$ follows directly by using Theorem \ref{T:InnerRate} with $k = k_{best}$ and $\gamma = 1/3$. The bound for $\epsilon_{\omega}$ follows by substituting $\gamma = 1/3$ into Theorem \ref{T:GkBound}.
    \end{proof}
    \begin{remark}[Rate preservation for the generalized FW gap] \label{R:GenFWGap}
        Extending the formula in Definition \ref{D:GenFWGap} to all $\bx \in \real^{n}$, if $X^{\ast}$ is compact, then Theorems \ref{T:GkBound} and \ref{T:SimulConvergenceRate} remain valid with the generalized FW gap ${\hat \GGG}_{OP}$ in place of the normalized one, with the constant multiplied by $R$ defined below. Indeed, since $\varphi$ is proper, lsc and convex with a nonempty compact set of minimizers, all its level sets are bounded \cite[Corollary 8.7.1]{R1970-B}. By Lemma \ref{L:FunctionBound}, the sequence $\left\{ \varphi\left(\bx^{k}\right) \right\}_{k \in \nn}$ is bounded, and thus so is $\Sequp{\bx}{k}$. Hence, $R := \sup_{\bx \in X^{\ast} , k \in \nn} \norm{\bx^{k} - \bx} < \infty$. Now, Proposition \ref{P:RegStepFW} with $\bx^{+} = \bx^{k + 1}$, $\bv = \bv^{k}$ and $\eta = \eta_{k}$, together with $\varphi\left(\bx^{k + 1}\right) \geq \varphi\left(\bx\right)$ for all $\bx \in X^{\ast}$, yields after taking the supremum over $\bx \in X^{\ast}$
        \begin{equation*}
            {\hat \GGG}_{OP}\left(\bx^{k + 1}\right) \leq \frac{2R}{\eta_{k}}\norm{G_{1/\beta}^{f_{k} , g_{k}}\left(\bv^{k}\right)}.
        \end{equation*}
        This is the bound on $\GGG_{OP}\left(\bx^{k + 1}\right)$ in Proposition \ref{P:RegStepFW}, multiplied by $R$. Therefore, the proofs of Theorems \ref{T:GkBound} and \ref{T:SimulConvergenceRate} apply verbatim. In particular, the rate $O(k^{-1/3})$ of Theorem \ref{T:SimulConvergenceRate} is preserved.
    \end{remark}    

\subsection{Gradient Mapping Rates under a H\"olderian Error Bound}
    The preceding result gives simultaneous rates for the inner objective gap and the normalized FW gap. Together, these quantities measure progress toward feasibility and outer stationarity in the sense of Definition \ref{D:EpsStat}. We next study the progress of the MABPG Algorithm in terms of the gradient mapping of the outer problem, which incorporates the constraint set $X^{\ast}$ into the proximal gradient step. The bound on the gradient mapping, which we proved in Proposition \ref{P:GradMapFWgap}, provides an additional interpretation of the accuracy of the generated iterates. However, in order to achieve the desired rate and since the iterates need not belong to $X^{\ast}$, we also need to control their distance to this set. A H\"olderian error bound on the inner objective provides this control. As we show below, it also improves the inner objective rate. 
    \begin{definition}[H\"olderian type error bound] \label{D:Holder}
        Let $p \geq 1$. A proper, lsc and convex function $h : \real^{n} \rightarrow \erl$ with $\argmin h \neq \emptyset$ satisfies a H\"olderian error bound of order $p$ with $\tau > 0$ if for all $\bx \in \real^{n}$
        \begin{equation*}
            \tau \cdot \dist\left(\bx , \argmin h\right)^{p} \leq h\left(\bx\right) - \min_{\bz \in \real^{n}} h\left(\bz\right).
        \end{equation*}
    \end{definition}
    The case $p = 2$ is the quadratic growth condition. Weak sharp minima ($p = 1$) hold, for example, for any polyhedral convex function with a nonempty set of minimizers \cite{BF1993}. Notice that the smaller $p$ is, the stronger the condition is near $\argmin h$.

    Throughout this subsection, we make the following additional assumption.
    \begin{assumption} \label{A:AssumpB}
        \begin{itemize}
            \item[$\rm{(i)}$] $\varphi$ satisfies the H\"olderian error bound of order $p \geq 1$ with parameter $\tau > 0$.
            \item[$\rm{(ii)}$] At least one of $\varphi$ and $\omega$ is coercive.
        \end{itemize}
    \end{assumption}
    For the rest of this section, we define $\bar{\bx}^{k} := P_{X^{\ast}}\left(\bx^{k}\right)$ for all $k \in \nn$.

    Now, we prove the following lemma.
    \begin{lemma} \label{L:SubGradBound}
        Let $\left\{ \left(\bx^{k} , \bv^{k} \right) \right\}_{k \in \nn}$ be a sequence generated by MABPG with $0 < \gamma \leq 2$. Then, the following items hold true.
        \begin{itemize}
            \item[$\rm{(i)}$] The sequences $\Sequp{\bx}{k}$ and $\Sequp{\bar{\bx}}{k}$ are bounded.
            \item[$\rm{(ii)}$] There exists $W > 0$ such that $\norm{\nabla\sigma\left(\bx^{k}\right) + \xi} \leq W$ for all $\xi \in \partial \psi\left(\bar{\bx}^{k}\right)$ and $k \in \nn$.
        \end{itemize}
    \end{lemma}
    \begin{proof}
        Since either $\varphi$ or $\omega$ is coercive (Assumption \ref{A:AssumpB}(ii)), the boundedness of $\Sequp{\bx}{k}$ follows from Lemma \ref{L:FunctionBound}, since both sequences $\left\{ \varphi\left(\bx^{k}\right) \right\}_{k \in \nn}$ and $\left\{ \omega\left(\bx^{k}\right) \right\}_{k \in \nn}$ are bounded. Moreover, since $P_{X^{\ast}}$ is nonexpansive, for any $\bx \in X^{\ast}$ we have $\norm{\bar{\bx}^{k} - \bx} = \norm{P_{X^{\ast}}\left(\bx^{k}\right) - P_{X^{\ast}}\left(\bx\right)} \leq \norm{\bx^{k} - \bx}$, and thus $\Sequp{\bar{\bx}}{k}$ is bounded as well. This proves the first item.
 
        Since $\nabla \sigma$ is continuous, it is bounded on bounded sets. Since $\psi$ is convex and real-valued (see Assumption \ref{A:AssumptionA}(iii)), its subdifferential is bounded on bounded sets (see \cite[Theorem 24.7]{R1970-B}). The second item now follows from the first one.
    \end{proof}
    The next lemma bounds the inner objective gap by the distance to $X^{\ast}$, and it is the common ingredient of all the results of this subsection.
    \begin{lemma} \label{L:InnerGapDist}
        Let $\left\{ \left(\bx^{k} , \bv^{k} \right) \right\}_{k \in \nn}$ be a sequence generated by MABPG with $0 < \gamma \leq 2$. Then, for all $k \geq 1$, we have
        \begin{equation} \label{L:InnerGapDist:1}
            \varphi\left(\bx^{k_{best} + 1}\right) - \varphi^{\ast} \leq \frac{2C + W}{k^{\alpha}} \cdot \dist\left(\bx^{k_{best} + 1} , X^{\ast}\right),
        \end{equation}
        where $\alpha := \min\left\{\gamma , (1 + \gamma)/2\right\}$, and $W$ as defined in Lemma \ref{L:SubGradBound}.
    \end{lemma}
    \begin{proof}
        For simplicity, we denote $s := k_{best}$, $G_{1/\beta} := G_{1/\beta}^{f_{s} , g_{s}}$ and $d := \norm{\bx^{s + 1} - \bar{\bx}^{s + 1}} = \dist\left(\bx^{s + 1} , X^{\ast}\right)$. Recall that $\bx^{s + 1} = T_{1/\beta}^{f_{s} , g_{s}}\left(\bv^{s}\right)$ (see Step \eqref{Bi-AGUpdate:3}) and that $\beta$ is a Lipschitz constant of $\nabla f_{s}$ (see Remark \ref{R:Lipfk}). Therefore, Proposition \ref{P:RegStepFW} with $\bx^{+} = \bx^{s + 1}$, $\bv = \bv^{s}$, $\bx = \bar{\bx}^{s + 1}$ and $\eta = \eta_{s}$, yields
        \begin{equation} \label{L:InnerGapDist:2}
            \varphi\left(\bx^{s + 1}\right) - \varphi\left(\bar{\bx}^{s + 1}\right) + \eta_{s}\left(\psi\left(\bx^{s + 1}\right) - \psi\left(\bar{\bx}^{s + 1}\right) + \act{\nabla \sigma\left(\bx^{s + 1}\right) , \bx^{s + 1} - \bar{\bx}^{s + 1}}\right) \leq 2\norm{G_{1/\beta}\left(\bv^{s}\right)} \cdot d.
        \end{equation}
        Let $\xi \in \partial \psi\left(\bar{\bx}^{s + 1}\right)$, which is nonempty since $\psi$ is real-valued. From the subgradient inequality, we have $\psi\left(\bx^{s + 1}\right) - \psi\left(\bar{\bx}^{s + 1}\right) \geq \act{\xi , \bx^{s + 1} - \bar{\bx}^{s + 1}}$. Thus, the Cauchy-Schwarz inequality and Lemma \ref{L:SubGradBound}(ii), yields
        \begin{equation*}
            \psi\left(\bx^{s + 1}\right) - \psi\left(\bar{\bx}^{s + 1}\right) + \act{\nabla \sigma\left(\bx^{s + 1}\right) , \bx^{s + 1} - \bar{\bx}^{s + 1}} \geq \act{\nabla \sigma\left(\bx^{s + 1}\right) + \xi , \bx^{s + 1} - \bar{\bx}^{s + 1}} \geq -W d.
        \end{equation*}
        Plugging this into \eqref{L:InnerGapDist:2} and recalling that $\varphi\left(\bar{\bx}^{s + 1}\right) = \varphi^{\ast}$, since $\bar{\bx}^{s + 1} \in X^{\ast}$, we obtain
        \begin{equation} \label{L:InnerGapDist:3}
            \varphi\left(\bx^{s + 1}\right) - \varphi^{\ast} \leq \left(2\norm{G_{1/\beta}\left(\bv^{s}\right)} + \eta_{s}W\right) \cdot d.
        \end{equation}
        It remains to bound the two terms in the parentheses. From Proposition \ref{P:GkBound}, we have
        \begin{equation*}
            \norm{G_{1/\beta}\left(\bv^{k_{best}}\right)} \leq \frac{C}{k^{\left(1 + \gamma\right)/2}} \leq \frac{C}{k^{\alpha}},
        \end{equation*}
        where the last inequality follows from the fact that $\alpha \leq (1 + \gamma)/2$ and $k \geq 1$. Moreover, for $s = k_{best} \geq k$, we have $\eta_{s} = \left(s + 1\right)^{-\gamma} \leq k^{-\gamma} \leq k^{-\alpha}$. Hence, $2\norm{G_{1/\beta}\left(\bv^{k_{best}}\right)} + \eta_{k_{best}}W \leq \left(2C + W\right)k^{-\alpha}$, and the result follows from \eqref{L:InnerGapDist:3}.
    \end{proof}
    \begin{remark} \label{R:Theta}
        For $\gamma \leq 1$ we have $\alpha = \gamma$. For $1 < \gamma \leq 2$ we have $\alpha = (1 + \gamma)/2 < \gamma$: the term $\eta_{s}W$ decays faster than the gradient mapping term, which then determines the rate.
    \end{remark}
    The H\"olderian error bound bounds the left-hand side of \eqref{L:InnerGapDist:1} from below by the distance to $X^{\ast}$. Thus, for $p > 1$, this yields rates for both the distance and the inner objective gap.
    \begin{theorem}[Improved inner rate] \label{T:ImprovedInnerRate}
        Let $\left\{ \left(\bx^{k} , \bv^{k} \right) \right\}_{k \in \nn}$ be a sequence generated by MABPG with $0 < \gamma \leq 2$, and suppose that Assumption \ref{A:AssumpB} holds with $p > 1$. Then, for any $\bx \in X^{\ast}$ and $k \geq 1$, we have
        \begin{equation*}
            \dist\left(\bx^{k_{best} + 1} , X^{\ast}\right) \leq \left(\frac{2C + W}{\tau}\right)^{\frac{1}{p - 1}} \frac{1}{k^{\alpha/(p - 1)}} \quad \text{and} \quad \varphi\left(\bx^{k_{best} + 1}\right) - \varphi\left(\bx\right) \leq \frac{\left(2C + W\right)^{p/(p - 1)}}{\tau^{1/(p - 1)}} \cdot \frac{1}{k^{\alpha p/(p - 1)}},
        \end{equation*}
        where $\alpha$ is as defined in Lemma \ref{L:InnerGapDist}.
    \end{theorem}
    \begin{proof}
        From Assumption \ref{A:AssumpB}(i) and Lemma \ref{L:InnerGapDist}, we have
        \begin{equation*}
            \tau \dist\left(\bx^{k_{best} + 1} , X^{\ast}\right)^{p} \leq \varphi\left(\bx^{k_{best} + 1}\right) - \varphi\left(\bx\right) \leq \frac{2C + W}{k^{\alpha}} \cdot \dist\left(\bx^{k_{best} + 1} , X^{\ast}\right).
        \end{equation*}
        If $\dist\left(\bx^{k_{best} + 1} , X^{\ast}\right) = 0$, then both bounds hold trivially. Otherwise, after dividing both sides by $\tau \dist\left(\bx^{k_{best} + 1} , X^{\ast}\right)$, we have that $\dist\left(\bx^{k_{best} + 1} , X^{\ast}\right)^{p - 1} \leq \left(2C + W\right)/\left(\tau k^{\alpha}\right)$, which is the first bound. Plugging it into the right-hand side above yields the second bound.
    \end{proof}
    For $p = 1$, the same argument shows that selected outputs of MABPG reach $X^{\ast}$ in finitely many iterations.
    \begin{theorem}[Finite feasibility under weak sharp minima] \label{T:FiniteFeasibility}
        Let $\left\{ \left(\bx^{k} , \bv^{k} \right) \right\}_{k \in \nn}$ be a sequence generated by MABPG with $0 < \gamma \leq 2$, and suppose that Assumption \ref{A:AssumpB} holds with $p = 1$. Let $\bx \in X^{\ast}$. Then, for all $k > \left(\left(2C + W\right)/\tau\right)^{1/\alpha}$, we have $\bx^{k_{best} + 1} \in X^{\ast}$, where $\alpha$ is as defined in Lemma \ref{L:InnerGapDist}.
    \end{theorem}
    \begin{proof}
        From Assumption \ref{A:AssumpB}(i) and Lemma \ref{L:InnerGapDist}, we have
        \begin{equation*}
            \tau\dist\left(\bx^{k_{best} + 1} , X^{\ast}\right) \leq \varphi\left(\bx^{k_{best} + 1}\right) - \varphi\left(\bx\right) \leq \frac{2C + W}{k^{\alpha}} \cdot \dist\left(\bx^{k_{best} + 1} , X^{\ast}\right).
        \end{equation*}
        Since $k > \left(\left(2C + W\right)/\tau\right)^{1/\alpha}$, we have $\left(2C + W\right)/k^{\alpha} < \tau$, and thus the inequality above can hold only if $\dist\left(\bx^{k_{best} + 1} , X^{\ast}\right) = 0$. Since $X^{\ast}$ is closed, this means that $\bx^{k_{best} + 1} \in X^{\ast}$.
    \end{proof}
    We are now ready to establish a rate for the gradient mapping of Problem \eqref{Prob:OP}. 
    \begin{theorem}[Gradient mapping rate] \label{T:GradientMappingRate}
        Let $\left\{ \left(\bx^{k} , \bv^{k}\right) \right\}_{k \in \nn}$ be a sequence generated by MABPG with $0 < \gamma < 1$. 
        \begin{itemize}
            \item[$\rm{(i)}$] If Assumption \ref{A:AssumpB} holds with $p > 1$, then there exists $D > 0$ such that for all $k \geq 1$, we have
            \begin{equation*}
                \norm{G_{1/\beta_{\sigma}}^{\sigma , \psi + \delta_{X^{\ast}}}\left(\bx^{k_{best} + 1}\right)} \leq \frac{D}{k^{\min\left\{1 - \gamma , \gamma/(p - 1)\right\}/2}}.
            \end{equation*}
            \item[$\rm{(ii)}$] If Assumption \ref{A:AssumpB} holds with $p = 1$, then for all $k > \left(\left(2C + W\right)/\tau\right)^{1/\gamma}$, we have
            \begin{equation*}
                \norm{G_{1/\beta_{\sigma}}^{\sigma , \psi + \delta_{X^{\ast}}}\left(\bx^{k_{best} + 1}\right)} \leq \frac{8C}{k^{(1 - \gamma)/2}}.
            \end{equation*}
        \end{itemize}
    \end{theorem}
    \begin{proof}
        For simplicity, we denote $G_{1/\beta_{\sigma}} := G_{1/\beta_{\sigma}}^{\sigma , \psi + \delta_{X^{\ast}}}$ and $d_{k} := \dist\left(\bx^{k_{best} + 1} , X^{\ast}\right)$. Applying Proposition \ref{P:GradMapFWgap} with $t = 1/\beta_{\sigma}$, $\bx = \bx^{k_{best} + 1}$, $\bx^{\ast}$ as its projection $\bar{\bx}^{k_{best} + 1}$ onto $X^{\ast}$, and $\xi \in \partial \psi\left(\bar{\bx}^{k_{best} + 1}\right)$, and then using Lemma \ref{L:SubGradBound}(ii), we obtain
        \begin{equation} \label{T:GradientMappingRate:1}
            \norm{G_{1/\beta_{\sigma}}\left(\bx^{k_{best} + 1}\right)} \leq 2\pos{\GGG_{OP}\left(\bx^{k_{best} + 1}\right)} + \beta_{\sigma} d_{k} + \sqrt{2\beta_{\sigma} W d_{k}}.
        \end{equation}
        From Theorem \ref{T:GkBound}, and the fact that $\left(k + 1\right)^{2 - \gamma} \geq 1$ (recall that $\gamma < 1$), we have
        \begin{equation} \label{T:GradientMappingRate:2}
            \pos{\GGG_{OP}\left(\bx^{k_{best} + 1}\right)} \leq \frac{4C}{k^{\left(1 - \gamma\right)/2}}.
        \end{equation}
        We first prove Item (i). From Theorem \ref{T:ImprovedInnerRate}, we have $d_{k} \leq \left(\left(2C + W\right)/\tau\right)^{1/(p - 1)}k^{-\gamma/(p - 1)}$. Plugging this bound and \eqref{T:GradientMappingRate:2} into \eqref{T:GradientMappingRate:1} yields
        \begin{align*}
            \norm{G_{1/\beta_{\sigma}}\left(\bx^{k_{best} + 1}\right)} & \leq \frac{8C}{k^{(1 - \gamma)/2}} + \beta_{\sigma}\left(\frac{2C + W}{\tau}\right)^{\frac{1}{p - 1}}\frac{1}{k^{\gamma/(p - 1)}} + \sqrt{2\beta_{\sigma} W}\left(\frac{2C + W}{\tau}\right)^{\frac{1}{2(p - 1)}}\frac{1}{k^{\gamma/(2(p - 1))}} \\
            & \leq \frac{D}{k^{\min\left\{1 - \gamma , \gamma/(p - 1)\right\}/2}},
        \end{align*}
        where the last inequality follows from the fact that each of the three powers of $k$ is bounded from above by $k^{-\min\left\{1 - \gamma , \gamma/(p - 1)\right\}/2}$ (recall that $k \geq 1$, and thus $k^{-\gamma/(p - 1)} \leq k^{-\gamma/(2(p - 1))}$).
\medskip
 
        Now we prove Item (ii). From Theorem \ref{T:FiniteFeasibility}, we have $d_{k} = 0$ for all $k > \left(\left(2C + W\right)/\tau\right)^{1/\gamma}$. Thus, the result follows from \eqref{T:GradientMappingRate:1} and \eqref{T:GradientMappingRate:2}.
    \end{proof}
    \begin{remark}[Choice of $\gamma$] \label{R:MABPGRates}
        \begin{itemize}
            \item[$\rm{(i)}$] For $p > 1$, the two exponents in Theorem \ref{T:GradientMappingRate}(i) are balanced by $\gamma = (p - 1)/p$. With this choice, $\alpha := \min\left\{\gamma , (1 + \gamma)/2\right\} = (p - 1)/p$, and thus Theorems \ref{T:ImprovedInnerRate} and \ref{T:GradientMappingRate}(i) give, at the same iterate $\bx^{k_{best} + 1}$, an inner objective gap of order $O(k^{-1})$ and an outer gradient mapping norm of order $O(k^{-1/(2p)})$. In particular, for $p = 2$ and $\gamma = 1/2$, these are $O(k^{-1})$ and $O(k^{-1/4})$.
            \item[$\rm{(ii)}$] For $p > 1$ and $0 < \gamma < 1$, the inner objective gap and the positive part of the normalized FW gap have exponents $p\gamma/\left(p - 1\right)$ and $\left(1 - \gamma\right)/2$, respectively. These exponents are balanced by $\gamma = \left(p - 1\right)/\left(3p - 1\right)$. With this choice, $\alpha := \min \left\{ \gamma , \left(1 + \gamma\right)/2 \right\} = \left(p - 1\right)/\left(3p - 1\right)$, and thus Theorems \ref{T:ImprovedInnerRate} and \ref{T:GkBound} give, at the same iterate $\bx^{k_{best} + 1}$, an inner objective gap and $\pos{\GGG_{OP}\left(\bx^{k_{best} + 1}\right)}$ of order $O(k^{-p/\left(3p - 1\right)})$. In particular, for $p = 2$ and $\gamma = 1/5$, both are of order $O(k^{-2/5})$.
            \item[$\rm{(iii)}$] For $p = 1$, the iterate $\bx^{k_{best} + 1}$ is feasible after finitely many iterations, and its outer gradient mapping norm is of order $O(k^{-(1 - \gamma)/2})$. Thus, small values of $\gamma$ give a rate arbitrarily close to $O(k^{-1/2})$, the rate of the proximal gradient method for single-level non-convex composite problems. The finite-feasibility threshold is $\left(\left(2C + W\right)/\tau\right)^{1/\gamma}$, which may be large for small values of $\gamma$.
            \item[$\rm{(iv)}$] In view of Theorem \ref{T:ImprovedInnerRate}, if only the inner problem matters, larger values of $\gamma$ are better. For $p = 2$, taking $\gamma = 2$ gives $\alpha = 3/2$ and thus $\dist\left(\bx^{k_{best} + 1} , X^{\ast}\right) = O(k^{-3/2})$ and $\varphi\left(\bx^{k_{best} + 1}\right) - \varphi^{\ast} = O(k^{-3})$. The price is the outer guarantee: Theorem \ref{T:GradientMappingRate} requires $\gamma < 1$, so for $\gamma \geq 1$ we have no rate for the outer gradient mapping.
        \end{itemize}
    \end{remark}    
    Remark \ref{R:MABPGRates} reveals a trade-off in the choice of $\gamma$. For $p = 2$, the choice $\gamma = 1/2$ gives the best outer gradient mapping rate established by Theorem \ref{T:GradientMappingRate}, namely $O(k^{-1/4})$, together with an inner objective gap of order $O(k^{-1})$. Increasing $\gamma$ improves the inner rate, but weakens the outer gradient mapping guarantee. In particular, $\gamma = 2$ gives an inner objective gap of order $O(k^{-3})$, but no outer gradient mapping rate is established for this choice. In the next subsection, we show that a projection phase after MABPG improves the inner objective gap from $O(k^{-1})$ to $O(k^{-3})$, while preserving the outer gradient mapping rate $O(k^{-1/4})$ at the same iterate.
    \begin{remark}[Computable stopping criterion] \label{R:StoppingCriterion}
        Evaluating the outer gradient mapping generally requires solving an auxiliary problem over $X^{\ast}$. In contrast, Step \eqref{Bi-AGUpdate:3} provides the residual
        \begin{equation*}
            r_{k} := \norm{G_{1/\beta}^{f_{k} , g_{k}}\left(\bv^{k}\right)} = \beta\norm{\bv^{k} - \bx^{k + 1}}
        \end{equation*}
        without additional oracle calls. Proposition \ref{P:RegStepFW} gives
        \begin{equation*}
            \pos{\GGG_{OP}\left(\bx^{k + 1}\right)} \leq \frac{2r_{k}}{\eta_{k}}.
        \end{equation*}
        Thus, if $\varphi^{\ast}$ is available, the conditions
        \begin{equation*}
            \Delta_{k} := \varphi\left(\bx^{k + 1}\right) - \varphi^{\ast} \leq \epsilon_{\varphi}, \qquad \frac{2r_{k}}{\eta_{k}} \leq \epsilon_{\omega}
        \end{equation*}
        certify that $\bx^{k + 1}$ is an $(\epsilon_{\varphi} , \epsilon_{\omega})$-stationary point as in Definition \ref{D:EpsStat}, without an error-bound assumption.

        Suppose, in addition, that Assumption \ref{A:AssumpB}(i) holds with known parameters $p , \tau$, and that $\psi$ is globally Lipschitz continuous with known constant $L_{\psi}$. Since every $\xi \in \partial\psi\left(P_{X^{\ast}}\left(\bx^{k + 1}\right)\right)$ satisfies $\norm{\xi} \leq L_{\psi}$, Proposition \ref{P:GradMapFWgap} yields, for every $t > 0$,
        \begin{equation} \label{E:ComputableBound}
            \norm{G_{t}^{\sigma , \psi + \delta_{X^{\ast}}}\left(\bx^{k + 1}\right)} \leq \frac{4r_{k}}{\eta_{k}} + \frac{1}{t}\left(\frac{\Delta_{k}}{\tau}\right)^{1/p} + \sqrt{\frac{2}{t}\left(\norm{\nabla\sigma\left(\bx^{k + 1}\right)} + L_{\psi}\right)}\left(\frac{\Delta_{k}}{\tau}\right)^{1/(2p)}.
        \end{equation}
        The right-hand side is computable without projecting onto $X^{\ast}$. In particular, when $\psi \equiv 0$, we take $L_{\psi} = 0$. If only $\ell \leq \varphi^{\ast}$ is available, both certificates remain valid with $\Delta_{k}$ replaced by $\varphi\left(\bx^{k + 1}\right) - \ell$.
    \end{remark}

\subsection{A Projection Phase for Feasibility} \label{SSec:MABPGP}
    Theorem \ref{T:ImprovedInnerRate} shows that, for $p = 2$ and $\gamma = 1/2$, the distance of $\bx^{k_{best} + 1}$ to $X^{\ast}$ is of order $O(k^{-1/2})$. In many applications, feasibility with respect to the inner problem is what matters most, and this can be improved substantially by projecting the output of MABPG onto $X^{\ast}$. The projection of a point $\bz \in \real^{n}$ onto $X^{\ast}$ is the solution of the simple bi-level problem
    \begin{equation} \label{P:ProjectionBiLevelProblem} \tag{P$_{\bz}$}
        \min_{\bx \in X^{\ast}} \left\{ h_{\bz}\left(\bx\right) := \frac{1}{2}\norm{\bx - \bz}^{2} \right\},
    \end{equation}
    whose outer objective is $1$-strongly convex with a $1$-Lipschitz continuous gradient. Thanks to the H\"olderian error bound (Assumption \ref{A:AssumpB}(i)), FBi-PG \cite{MST2026} solves this problem at fast rates. This leads to the following two-phase scheme. After terminating the first phase, its selected output is used to define the projection problem. The second-phase budget is then chosen in proportion to the number of first-phase iterations. Thus, the projection phase is a separate postprocessing step, while the updates of MABPG remain independent of the stopping horizon. 
\medskip
    

    \noindent\fbox{\parbox{\textwidth}{
    {\bf MABPG-Projected (MABPG-P)}
    \begin{itemize}[leftmargin=*]
        \item[] {\bf Input:} $k \in \nn$, $\gamma = 1/2$, $\beta = \beta_{f} + \beta_{\sigma}$, and $2 \leq a \in \nn$.
        \item[] {\bf Phase 1:} Run MABPG on Problem \eqref{Prob:OP} for $2k + 1$ iterations, and set $\bz = \bx^{k_{best} + 1}$.
        \item[] {\bf Phase 2:} Run FBi-PG \cite{MST2026} on Problem \eqref{P:ProjectionBiLevelProblem} for $2k$ iterations, with parameters $a$ and $\mu = 3/2$, step size $1/\tilde{\beta}$ where $\tilde{\beta} := \max\left\{ \beta , \beta_{f} + 1 \right\}$, and starting point $\bw^{0} = \bz$.
    \end{itemize}}}
\vspace{0.15in}

    In Phase 2, FBi-PG is applied with the regularization parameters $\alpha_{s} = \left(s + a\right)^{-3/2}$ and the momentum parameters $\theta_{s} = \left(s + a\right)/a$, $s \geq 0$. More precisely, starting from $\bw^{-1} = \bw^{0} = \bz$, and setting $\theta_{-1} = 0$, the sequence $\Sequp{\bw}{s}$ is generated, for any $s \geq 0$, by
    \begin{align*}
        \by^{s} & = \bw^{s} + \frac{\theta_{s - 1} - 1}{\theta_{s}}\left(\bw^{s} - \bw^{s - 1}\right), \\
        \bw^{s + 1} & = T_{1/\tilde{\beta}}^{f + \alpha_{s}h_{\bz} , g}    \left(\by^{s}\right).
    \end{align*}
    The second phase improves the inner rate from $O(k^{-1})$ to $O(k^{-3})$, while the rate of the gradient mapping is preserved. The proof is given in Appendix \ref{App:ProofTotalRate}.
    \begin{theorem}[Rates of MABPG-P] \label{T:TotalRate}
        Suppose that Assumption \ref{A:AssumpB} holds with $p = 2$. Let $\left\{ \left(\bx^{k} , \bv^{k}\right) \right\}_{k \in \nn}$ and $\Sequp{\bw}{s}$ be the sequences generated by Phase 1 and Phase 2 of MABPG-P, respectively. Then, there exist $\Theta_{1} , \Theta_{2} > 0$ such that for all $\bx \in X^{\ast}$ and $k \geq 1$, we have
        \begin{equation*}
            \varphi\left(\bw^{2k}\right) - \varphi\left(\bx\right) \leq \frac{\Theta_{1}}{k^{3}} \quad \text{and} \quad \norm{G_{1/\beta_{\sigma}}^{\sigma , \psi + \delta_{X^{\ast}}}\left(\bw^{2k}\right)} \leq \frac{\Theta_{2}}{k^{1/4}}.
        \end{equation*}
        Moreover, $\dist\left(\bw^{2k} , X^{\ast}\right) = O\left(k^{-3/2}\right)$.
    \end{theorem}
    
\section{Numerical Experiments} \label{Sec:Numerics}
    We evaluate MABPG and MABPG-P on simple bi-level problems with nonconvex outer objectives. The experiments examine their performance in reducing the inner objective gap and outer gradient mapping norm, and the effect of the projection phase in MABPG-P. We compare our methods with CG-BiO \cite{JAMH2023} and $\text{IPR-VFISTA}_{\text{nc}}$ \cite{SBY2025}. We also include FISTA \cite{BT2009} applied only to the inner problem as a reference for the performance obtained without using the outer objective.
\medskip

    In each of the two experiments below, all methods start from the same randomly generated point within each run and are given an optimization budget of $T = 1000$ seconds. We report averages over five runs. MABPG uses $\eta_{k} = \left(k + 1\right)^{-\gamma}$, step size $1/\beta$ with $\beta = \beta_{f} + \beta_{\sigma}$, and the additional candidate $\bu^{k + 1}$ described in Remark \ref{R:AddCand}; We record the proximal-gradient candidates $\bx^{k + 1}$ analyzed in Section \ref{Sec:MABPG}. For MABPG-P, we run MABPG for the first $500$ seconds, pass its final proximal-gradient iterate $\bx^{k + 1}$ as the reference point $\bz$ for the projection phase \eqref{P:ProjectionBiLevelProblem}, and then run FBi-PG \cite{MST2026} for the remaining $500$ seconds (recording $\bw^{s}$) with $\alpha_{s} = \left(s + a\right)^{-3/2}$, $\theta_{s} = \left(s + a\right)/a$, $a = 2$, and step size $1/\tilde{\beta}$, where $\tilde{\beta} = \max\left\{\beta_{f} + \beta_{\sigma} , \beta_{f} + 1\right\}$. Thus, compared to Theorem \ref{T:TotalRate}, the experimental implementation of MABPG-P differs in three respects: it uses $\gamma = 4/5$ (rather than $\gamma = 1/2$), a final-iterate handoff $\bz = \bx^{k + 1}$ (rather than $\bx^{k_{best} + 1}$), and an equal $500$-second time split between phases (rather than prescribed iteration counts $2k + 1$ and $2k$). For each experiment and each baseline algorithm, the main comparison displays the configuration with which that algorithm performed best on that experiment in terms of final outer gradient mapping norm $\norm{G_{t}\left(\bx^{k}\right)}$, while MABPG and MABPG-P use the regularization exponent $\gamma$ that achieved the best performance for vanilla MABPG on that experiment ($\gamma = 4/5$ in both experiments); all tested hyperparameter comparisons across both experiments are reported in Appendix \ref{app:hyperparameter_configs}. All benchmark experiments are implemented in Python using NumPy and SciPy and executed on Linux compute nodes of the Technion Zeus HPC cluster, with each algorithm restricted to a single CPU thread ($16\,\mathrm{GB}$ RAM allocated per job). To ensure consistent timing across methods, the $1000$-second budget counts optimization wall-clock time only (including the preliminary Frank-Wolfe phase of CG-BiO, which is part of that method); one-time problem precomputations (Lipschitz constants and matrix factorizations) are performed beforehand, and the timer is paused during all diagnostic objective and gradient-mapping evaluations.
\medskip

    Since $\psi \equiv 0$ in both experiments, we evaluate inner accuracy through the objective gap $\varphi\left(\bx^{k}\right) - \varphi^{\ast}$ and outer stationarity through the gradient mapping norm with respect to $\sigma$ and $\delta_{X^{\ast}}$, shortly written here as $\norm{G_{t}\left(\bx^{k}\right)}$, with $t = 1/\beta_{\sigma}$ (matching Theorems \ref{T:GradientMappingRate} and \ref{T:TotalRate}), where $\bx^{k}$ denotes the recorded iterate of a given method (and $\bw^{s}$ during Phase 2 of MABPG-P).
\medskip

    Both experiments use a smooth outer objective, with $\psi \equiv 0$ and $\sigma$ given by the Cauchy regularizer
    \begin{equation*}
        \sigma\left(\bx\right) := \sum_{i = 1}^{d} \frac{\lambda\rho^{2}}{2}\log\left(1 + \frac{x_{i}^{2}}{\rho^{2}}\right).
    \end{equation*}
    This penalty behaves quadratically near zero, while its derivative decays to zero as the magnitude of a coordinate increases. It therefore provides a smooth nonconvex criterion for selecting among inner solutions. We set $\lambda = 0.1$ and $\rho = 1$, for which the exact gradient Lipschitz constant is $\beta_{\sigma} = \lambda = 0.1$.
    
\subsection{Least Squares with a Nonconvex Outer Regularizer} \label{Subsec:LeastSquares}
    We consider the simple bi-level problem
    \begin{equation*}
        \min_{\bx \in X^{\ast}} \sigma\left(\bx\right), \qquad X^{\ast} := \argmin_{\bx \in \mathcal{Z}} \left\{ f\left(\bx\right) := \frac{1}{2}\norm{A\bx - b}^{2} \right\},
    \end{equation*}
    where $\sigma$ is the Cauchy regularizer defined above and $\mathcal{Z} := \left\{ \bx \in \real^{n} : \norm{\bx} \leq 2 \right\}$. Thus, the inner objective is $\varphi = f + \delta_{\mathcal{Z}}$.
\medskip

    We use the pre-constructed feature matrix and response vector from the Wikipedia Math Essential dataset \cite{RSHP2021} ($1068$ instances and $730$ features) without additional rescaling, and sample $400$ rows uniformly without replacement using a fixed random seed ($42$) to construct $A$ and $b$. The same sampled instance $\left(A , b\right)$ is used across all five runs, while each run draws a fresh random starting point $\bx^{0}$ with $\norm{\bx^{0}} = 1$ shared across all methods. The exact inner gradient Lipschitz constant $\beta_{f} = \norm{A}_{2}^{2} \approx 4007.8$ is computed from the largest singular value of $A$. The common algorithmic settings and computational budgets are described above.
\medskip

    Let $\bx^{\dagger} := A^{\dagger}b$, where $A^{\dagger}$ denotes the Moore-Penrose pseudoinverse. In our constructed instance, $A$ has full row rank (so $A\bx^{\dagger} = b$) and $\norm{\bx^{\dagger}} \approx 0.30 < 2$; consequently, $\varphi^{\ast} = 0$ and
    \begin{equation*}
        X^{\ast} = \left\{ \bx^{\dagger} + \bv : \bv \in \ker A ,\, \norm{\bx^{\dagger} + \bv} \leq 2 \right\}.
    \end{equation*}
    Since $A$ has more columns than rows, its nullspace is nontrivial, so $X^{\ast}$ contains infinitely many points and the outer objective provides an additional criterion for selecting among exact least-squares solutions.
\medskip

    To evaluate outer stationarity, we compute the projection of $\bu := \bx - t\nabla\sigma\left(\bx\right)$ onto $X^{\ast}$, with $t = 1/\beta_{\sigma}$. Since $A\bx^{\dagger} = b$, the affine set $\left\{ \by : A\by = b \right\}$ is nonempty, and the projection of $\bu$ onto it is
    \begin{equation*}
        \by = \bu + A^{\dagger}\left(b - A\bu\right).
    \end{equation*}
    Whenever $\norm{\by} \leq 2$, this point also equals $P_{X^{\ast}}\left(\bu\right)$, and the gradient mapping norm is $\norm{\bx - \by}/t$. We check this condition at every evaluation point; in the reported runs, it holds for all evaluated methods. 
    
\paragraph{Choice of the regularization exponent.}
    We first compare vanilla MABPG with $\gamma \in \left\{1/2 , 2/3 , 4/5\right\}$, keeping all other settings fixed (see Figure \ref{fig:gamma_comparison} in Appendix \ref{app:hyperparameter_configs}). The choice $\gamma = 1/2$ balances the exponents in the outer gradient mapping bound under quadratic growth, whereas the larger values give faster decay of the regularization parameter. On vanilla MABPG, $\gamma = 4/5$ attains the smallest final outer gradient mapping norm (Figure \ref{fig:gamma_comparison}); this confirms that the exponent optimizing the asymptotic bound need not give the best finite-time performance. We therefore use $\gamma = 4/5$ for both MABPG and its projected extension MABPG-P in the main comparison below.

\paragraph{Comparison with existing methods.}
    We compare MABPG and MABPG-P against FISTA and the best-performing configurations of CG-BiO ($c = 0.005$) and $\text{IPR-VFISTA}_{\text{nc}}$ ($\bar{\eta} = 1$ with outer step size $h/10$), selected from the parameter sweeps reported in Appendix \ref{app:hyperparameter_configs}.
\medskip

\begin{figure}[htbp]
        \centering
        \begin{subfigure}[t]{0.48\textwidth}
            \centering
            \includegraphics[width=\textwidth]{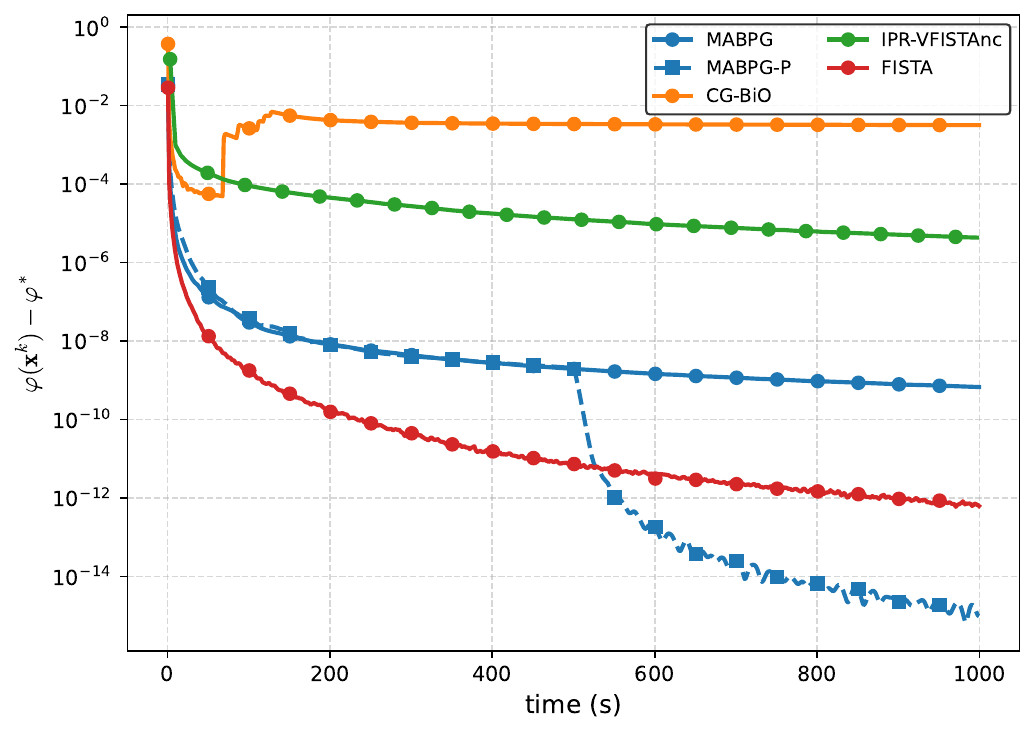}
            \caption{Inner-level objective gap $\varphi\left(\bx^{k}\right) - \varphi^{\ast}$.}
            \label{fig:g_gap_1}
        \end{subfigure}
        \hfill
        \begin{subfigure}[t]{0.48\textwidth}
            \centering
            \includegraphics[width=\textwidth]{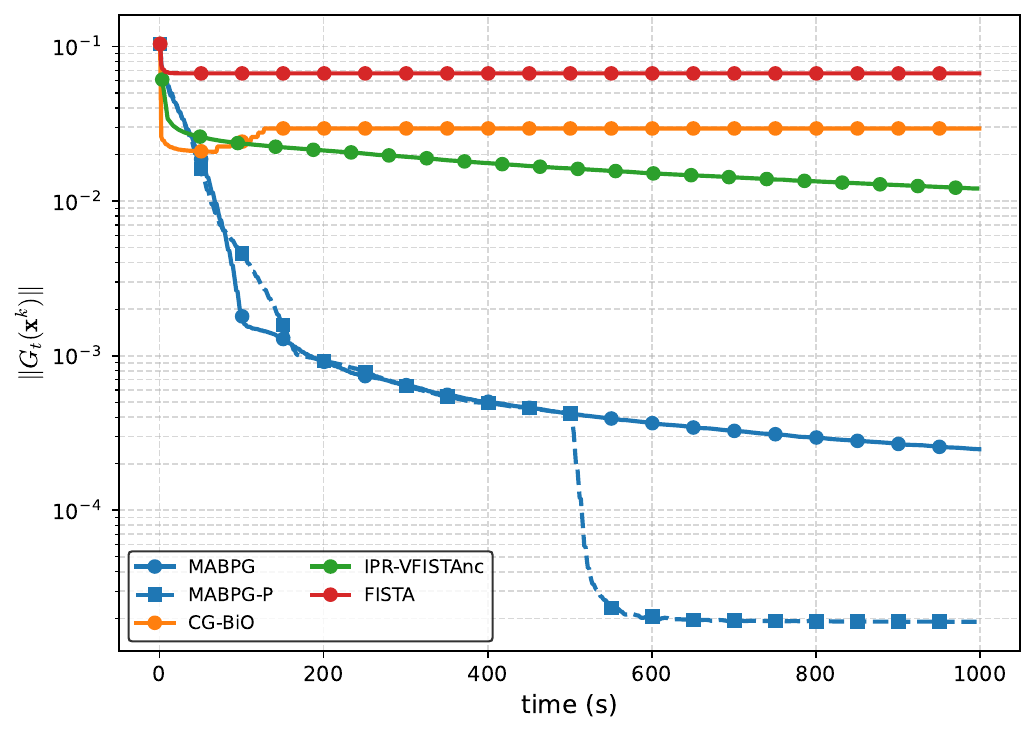}
            \caption{Norm of the gradient mapping $\norm{G_{t}\left(\bx^{k}\right)}$.}
            \label{fig:grad_mapping_1}
        \end{subfigure}
        \caption{Performance comparison of MABPG and MABPG-P against the best-performing configurations of baselines on the non-convex regularized least squares problem. Lines represent the average across 5 independent runs.}
        \label{fig:experiment_1}
    \end{figure}

    Figure \ref{fig:experiment_1} reports the inner objective gap and outer gradient mapping norm against elapsed time, averaged over five runs. FISTA initially reduces the inner objective gap faster than the bi-level methods, but its outer gradient mapping norm quickly plateaus because it ignores the outer objective. In contrast, MABPG steadily reduces both quantities and outperforms both CG-BiO and $\text{IPR-VFISTA}_{\text{nc}}$ in both metrics within the same time budget.
\medskip

    The projection phase further improves MABPG-P. After the switch at $500$ seconds, the projection phase sharply reduces both the inner objective gap and the outer gradient mapping norm, yielding better inner accuracy and outer stationarity than continuing MABPG for the full budget.
\medskip

    We also examine the empirical tightness of the computable stopping bound \eqref{E:ComputableBound} in Remark \ref{R:StoppingCriterion}. For this least-squares problem, we use $p = 2$, $\varphi^{\ast} = 0$, and $t = 1/\beta_{\sigma} = 10$. Whenever the affine projection $P_{\left\{\by : A\by = b\right\}}\left(\bx\right)$ belongs to $\mathcal{Z}$, it coincides with $P_{X^{\ast}}\left(\bx\right)$, so $\dist\left(\bx , X^{\ast}\right) = \dist\left(\bx , \left\{\by : A\by = b\right\}\right)$ and therefore the quadratic-growth inequality $\varphi\left(\bx\right) - \varphi^{\ast} = \frac{1}{2}\norm{A\bx - b}^{2} \geq \tau\,\dist\left(\bx , X^{\ast}\right)^{2}$ holds with $\tau = \sigma_{\min,+}^{2}\left(A\right)/2 \approx 2.51 \times 10^{-5}$, where $\sigma_{\min,+}\left(A\right)$ is the smallest positive singular value of $A$. We verify that this condition holds at all evaluated iterates. The ratio of the largest to smallest positive eigenvalues of $A^{\top}A$ is approximately $7.99 \times 10^{7}$.
    \medskip

\begin{figure}[htbp]
        \centering
        \begin{subfigure}[t]{0.48\textwidth}
            \centering
            \includegraphics[width=\textwidth]{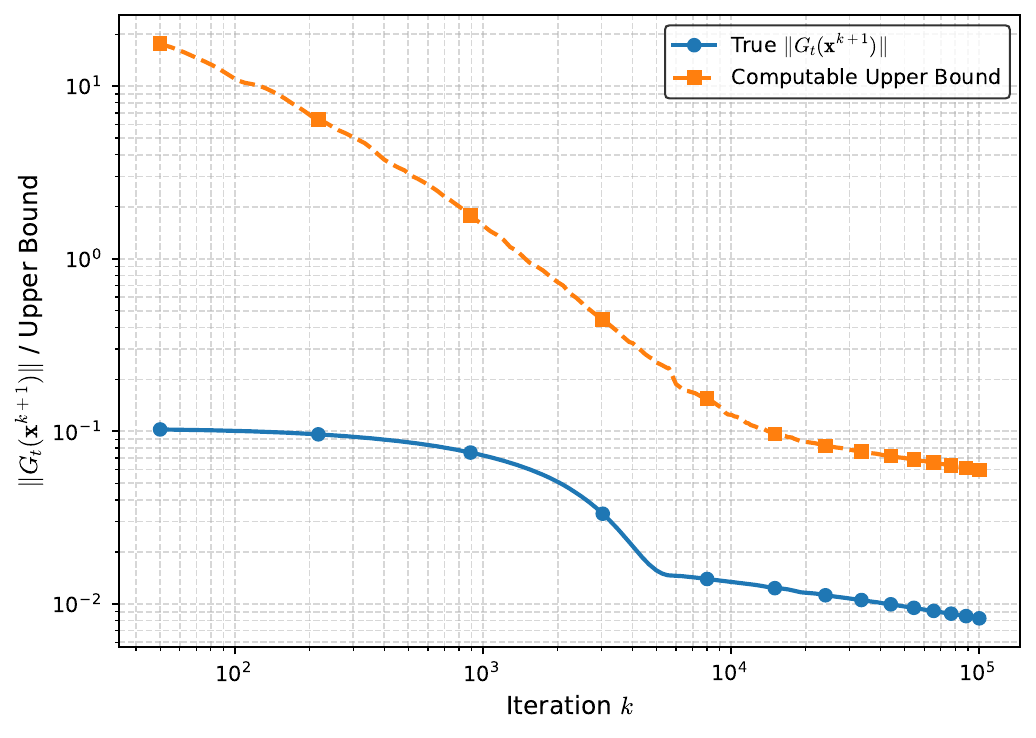}
            \caption{Versus iteration count $k$.}
            \label{fig:stopping_bound_curves}
        \end{subfigure}
        \hfill
        \begin{subfigure}[t]{0.48\textwidth}
            \centering
            \includegraphics[width=\textwidth]{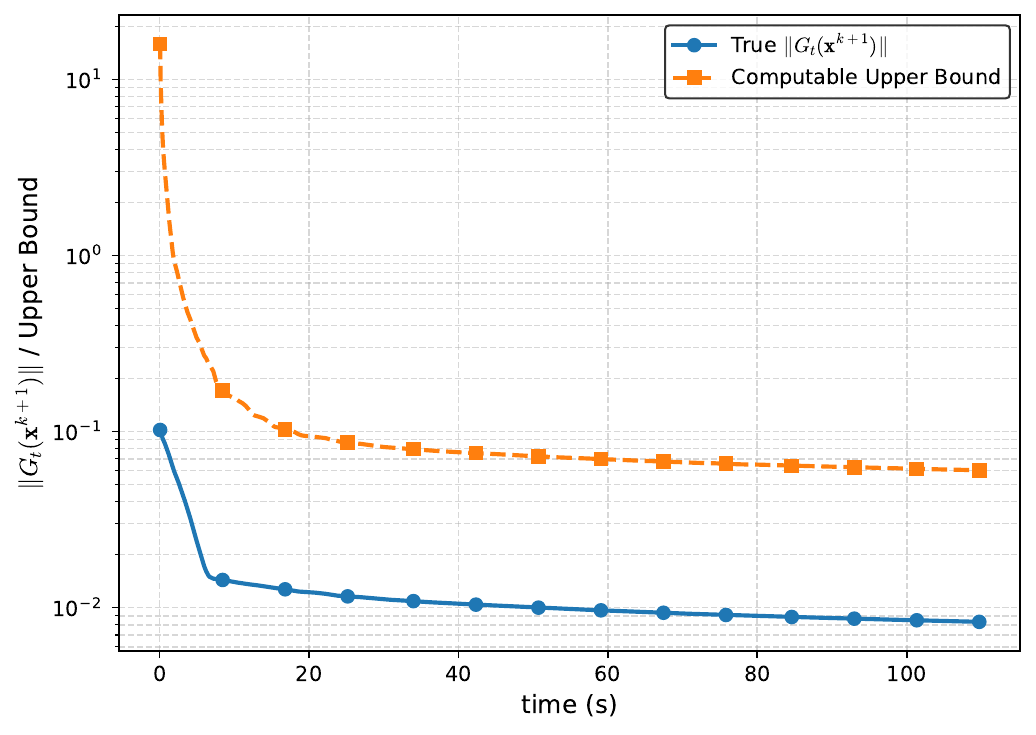}
            \caption{Versus elapsed time (seconds).}
            \label{fig:stopping_bound_time}
        \end{subfigure}
        \caption{True $\norm{G_{t}\left(\bx^{k + 1}\right)}$ versus the computable upper bound \eqref{E:ComputableBound} for MABPG with $\gamma = 1/2$ on the non-convex regularized least squares problem (averaged over 5 runs).}
        \label{fig:stopping_bound}
    \end{figure}

    Figure \ref{fig:stopping_bound} compares the computable bound with the outer gradient mapping norm versus both iteration count and elapsed time over the first $10^{5}$ iterations of MABPG with $\gamma = 1/2$ (the theoretically optimal exponent for the outer gradient mapping bound under quadratic growth; see Remark \ref{R:MABPGRates}(i)), averaged over five runs executed single-threaded in Python 3.12 with NumPy 2.4 on a separate Windows 11 workstation (Intel Core i7-1195G7 CPU @ 2.90\,GHz, $16\,\mathrm{GB}$ RAM); accordingly, the elapsed times in Figure \ref{fig:stopping_bound_time} are not intended to be compared directly with the cluster wall-clock times in Figures \ref{fig:experiment_1} and \ref{fig:experiment_2}. The ratio between the bound and the gradient mapping norm decreases from $109.8$ at $k = 100$ to $21.6$ at $k = 1000$ and reaches approximately $7.25$ at $k = 10^{5}$. At the final iteration, the bound and the gradient mapping norm are $5.98 \times 10^{-2}$ and $8.24 \times 10^{-3}$, respectively. Thus, despite the ill-conditioning, the bound becomes substantially more informative over the run and is within one order of magnitude of the gradient mapping norm at its end.

\subsection{Over-parameterized Text Classification} \label{Subsec:TextClass}
    We next consider an over-parameterized sparse text classification problem with the squared-hinge loss and the same Cauchy outer regularizer ($\lambda = 0.1$, $\rho = 1$, $\beta_{\sigma} = 0.1$):
    \begin{equation} \label{eq:sqhinge_bilevel}
        \min_{\bx \in X^{\ast}} \sigma\left(\bx\right), \qquad X^{\ast} := \argmin_{\bx \in \mathcal{Z}} \left\{ f\left(\bx\right) := \frac{1}{2n}\sum_{i = 1}^{n}\max\left\{0 ,\, 1 - b_{i}A_{i}^{\top}\bx\right\}^{2} \right\},
    \end{equation}
    where $A_{i}^{\top}$ is the $i$-th row of the feature matrix $A \in \real^{n \times d}$, $b_{i} \in \left\{-1 , 1\right\}$, and $\mathcal{Z} := \left\{ \bx \in \real^{d} : \norm{\bx} \leq 50 \right\}$, so that $\varphi := f + \delta_{\mathcal{Z}}$.
\medskip

    We construct the instance from the official LYRL2004 training split of the RCV1-v2 text categorization benchmark \cite{LYRL2004} (unit-norm TF-IDF feature vectors). Documents carrying at least one category in $\left\{\texttt{CCAT}, \texttt{ECAT}\right\}$ and none in $\left\{\texttt{GCAT}, \texttt{MCAT}\right\}$ are labeled $+1$; those carrying at least one category in $\left\{\texttt{GCAT}, \texttt{MCAT}\right\}$ and none in $\left\{\texttt{CCAT}, \texttt{ECAT}\right\}$ are labeled $-1$; documents belonging to both groups or neither are discarded. After deduplicating identical sparse feature vectors, we sample $n = 500$ documents uniformly without replacement using a fixed seed ($42$) and retain the $d = 7666$ nonzero feature columns ($36{,}494$ nonzero entries); the same sampled instance is reused across all five runs, and the inner gradient Lipschitz constant $\beta_{f} = \norm{A}_{2}^{2}/n \approx 0.0232$ is computed via sparse singular value decomposition. Since the rows of $A$ are linearly independent and the minimum-norm separating classifier satisfies $\norm{\bx_{\text{svm}}} \approx 16.60 < 50$, the exact inner optimum is $\varphi^{\ast} = 0$ and the lower-level solution set $X^{\ast} = \left\{ \bx \in \mathcal{Z} : b_{i}A_{i}^{\top}\bx \geq 1 \text{ for all } i = 1,\dots,n \right\}$ is a full-dimensional convex set with nonempty interior containing infinitely many separating classifiers. The Cauchy outer regularizer therefore selects among zero-loss classifiers.
\medskip

    To evaluate outer stationarity, we compute the projection of $\bu := \bx - t\nabla\sigma\left(\bx\right)$ onto $X^{\ast}$, with $t = 1/\beta_{\sigma} = 10$. Let $M \in \real^{n \times d}$ have $i$-th row $b_{i}A_{i}^{\top}$. Computing the projection $\by$ of $\bu$ onto the polyhedron $\left\{ \by \in \real^{d} : M\by \geq \mathbf{1} \right\}$ is a standard least-distance problem; since $M$ has full row rank $n = 500$, it reduces via a precomputed Cholesky factorization of $MM^{\top}$ to an $n$-dimensional nonnegative least-squares problem \cite[Chapter 23]{LH1995}, solved via SciPy's active-set Lawson--Hanson routine (\texttt{scipy.optimize.nnls}). Whenever $\norm{\by} \leq 50$, this point also equals $P_{X^{\ast}}\left(\bu\right)$, and the gradient mapping norm is $\norm{\bx - \by}/t$. We check this condition at every evaluation point; in the reported runs, it holds for all evaluated methods.
\medskip

    Within each run, all methods start from the same randomly generated point $\bx^{0}$ with $\norm{\bx^{0}} = 1$ and are given $T = 1000$ seconds. Following the parameter comparisons in Appendix \ref{app:hyperparameter_configs}, we use $\gamma = 4/5$ for MABPG and its projected extension MABPG-P (which achieves the smallest final $\norm{G_{t}\left(\bx^{k}\right)}$ for vanilla MABPG here), $c = 0.005$ for CG-BiO, and $\bar{\eta} = 0.01$ with outer step size $h = 1/\left(2\beta_{\sigma}\right)$ for $\text{IPR-VFISTA}_{\text{nc}}$.

\begin{figure}[htbp]
        \centering
        \begin{subfigure}[t]{0.48\textwidth}
            \centering
            \includegraphics[width=\textwidth]{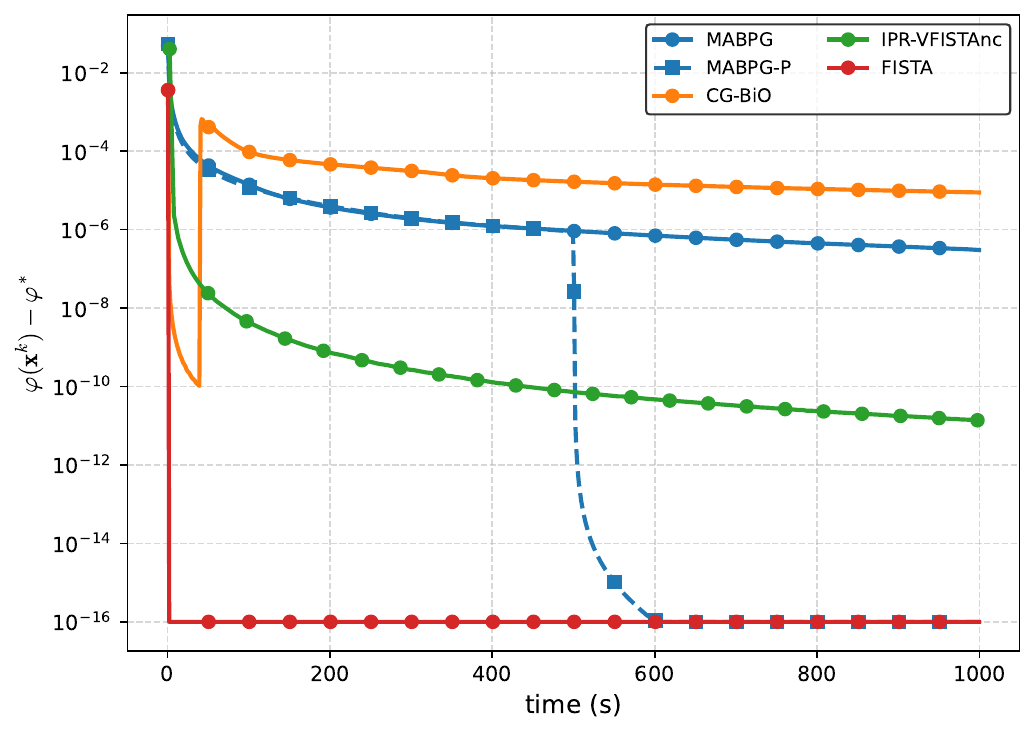}
            \caption{Inner-level objective gap $\varphi\left(\bx^{k}\right) - \varphi^{\ast}$.}
            \label{fig:g_gap_2}
        \end{subfigure}
        \hfill
        \begin{subfigure}[t]{0.48\textwidth}
            \centering
            \includegraphics[width=\textwidth]{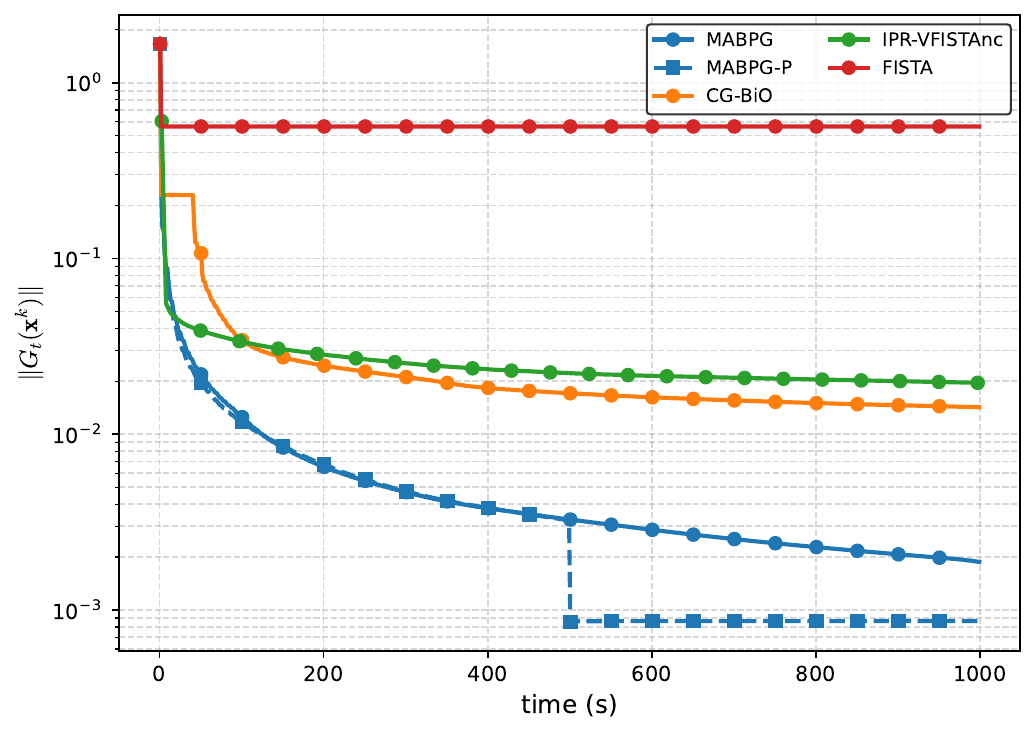}
            \caption{Norm of the gradient mapping $\norm{G_{t}\left(\bx^{k}\right)}$.}
            \label{fig:grad_mapping_2}
        \end{subfigure}
        \caption{Performance comparison of MABPG and MABPG-P against the best-performing configurations of baselines on the over-parameterized RCV1-v2 text classification problem. Lines represent the average across 5 independent runs.}
        \label{fig:experiment_2}
    \end{figure}

\paragraph{Comparison results.}
    Figure \ref{fig:experiment_2} shows both the inner objective gap $\varphi\left(\bx^{k}\right) - \varphi^{\ast}$ and the outer gradient mapping norm $\norm{G_{t}\left(\bx^{k}\right)}$ against elapsed time, averaged over five runs. FISTA reaches $\varphi\left(\bx^{k}\right) = 0$ almost immediately (within a fraction of a second, floored at $10^{-16}$ on the logarithmic plot) as its iterates satisfy all margin constraints, but plateaus at a large outer gradient mapping norm because it ignores the outer regularizer. Both CG-BiO and $\text{IPR-VFISTA}_{\text{nc}}$ improve substantially over FISTA, with $\text{IPR-VFISTA}_{\text{nc}}$ achieving a smaller inner gap and CG-BiO attaining a slightly smaller outer gradient mapping norm. MABPG outperforms both baselines in outer gradient mapping norm by roughly an order of magnitude, and the Phase 2 projection in MABPG-P drives the inner objective gap to machine precision while further reducing $\norm{G_{t}\left(\bx^{k}\right)}$ below that of continuing MABPG.

\appendix

\section{Hyperparameter Comparisons} \label{app:hyperparameter_configs}
    Because $T_{t}\left(\bx\right) \in X^{\ast}$ gives $\dist\left(\bx , X^{\ast}\right) \leq t\norm{G_{t}\left(\bx\right)}$, the outer gradient mapping norm $\norm{G_{t}\left(\bx^{k}\right)}$ jointly captures inner feasibility and outer stationarity and serves as our selection criterion across both benchmarks (inner gaps for the selected methods are shown in Figures \ref{fig:experiment_1} and \ref{fig:experiment_2}). All configurations use the common initializations and computational budgets of Section \ref{Sec:Numerics}, and curves represent averages over five runs.

    Across both experiments, we test the following parameter settings:
    \begin{itemize}[leftmargin=*]
        \item \textbf{MABPG:} $\eta_{k} = \left(k + 1\right)^{-\gamma}$ with $\gamma \in \left\{1/2 , 2/3 , 4/5\right\}$ (MABPG-P inherits the exponent selected for vanilla MABPG).
        \item \textbf{CG-BiO \cite{JAMH2023}:} analytical LMO over the intersection of the Euclidean ball $\mathcal{Z}$ and the cutting-plane halfspace, with step-size schedule $\gamma_{k} = \min\left\{1 , c/\sqrt{k + 1}\right\}$ for $c \in \left\{0.005 , 0.01 , 0.1 , 1\right\}$ (labeled $\gamma$ in the legends). The preliminary Frank-Wolfe phase uses step sizes $2/\left(k + 2\right)$ and stops when the inner Frank-Wolfe gap falls below $10^{-5}$ or when a cap of $10^{5}$ iterations is reached; if the iteration cap is reached before the gap tolerance (as occurs on the RCV1-v2 instance), the $10^{5}$-th iterate is passed as the starting point for the main cutting-plane phase.
        \item \textbf{$\text{IPR-VFISTA}_{\text{nc}}$ \cite{SBY2025}:} update scheme of Algorithm 3 in \cite{SBY2025}, using tuned constant outer step sizes rather than a horizon-dependent step size for the wall-clock comparison. At outer iteration $k$, the inner loop runs for $J_{k} = \max\left\{2 , \left\lceil\left(k + 1\right)^{2.01}\right\rceil\right\}$ iterations with regularization parameter $\eta_{k} = 16\left(\beta_{f} + \bar{\eta}\right)\left(\ln J_{k} / J_{k}\right)^{2}$, inner step size $1/\left(\beta_{f} + \eta_{k}\right)$, and strongly convex FISTA momentum $\left(\sqrt{\kappa_{k}} - 1\right)/\left(\sqrt{\kappa_{k}} + 1\right)$ where $\kappa_{k} = \left(\beta_{f} + \eta_{k}\right)/\eta_{k}$. Setting $h := 1/\left(2\beta_{\sigma}\right)$ (denoted $\hat{\gamma}$ in the legends), we test $\bar{\eta} \in \left\{0.01 , 0.1 , 1\right\}$ with constant outer step size $h$, as well as $h/2$, $h/10$, and $h/100$ with $\bar{\eta} = 1$.
    \end{itemize}

\begin{figure}[!htbp]
        \centering
        \begin{subfigure}[t]{0.48\textwidth}
            \centering
            \includegraphics[width=\textwidth]{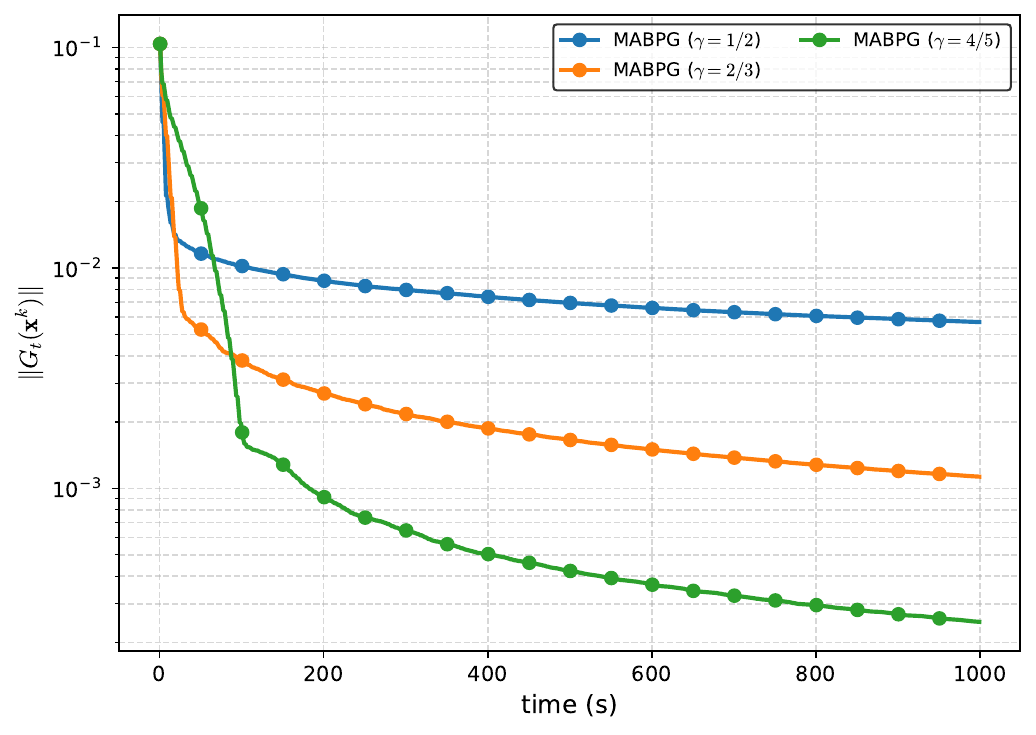}
            \caption{MABPG with $\gamma \in \{1/2, 2/3, 4/5\}$.}
            \label{fig:gamma_comparison}
        \end{subfigure}
        \hfill
        \begin{subfigure}[t]{0.48\textwidth}
            \centering
            \includegraphics[width=\textwidth]{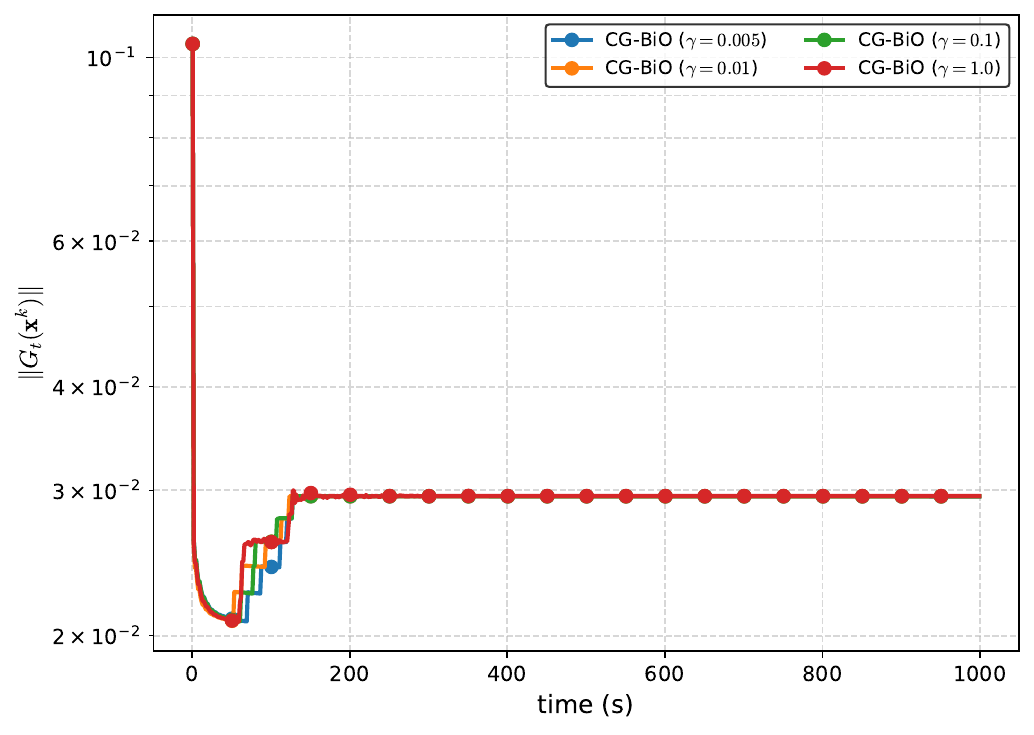}
            \caption{CG-BiO step-size coefficients.}
            \label{fig:app_least_squares_cgbio}
        \end{subfigure}
        \\[0.8em]
        \begin{subfigure}[t]{0.48\textwidth}
            \centering
            \includegraphics[width=\textwidth]{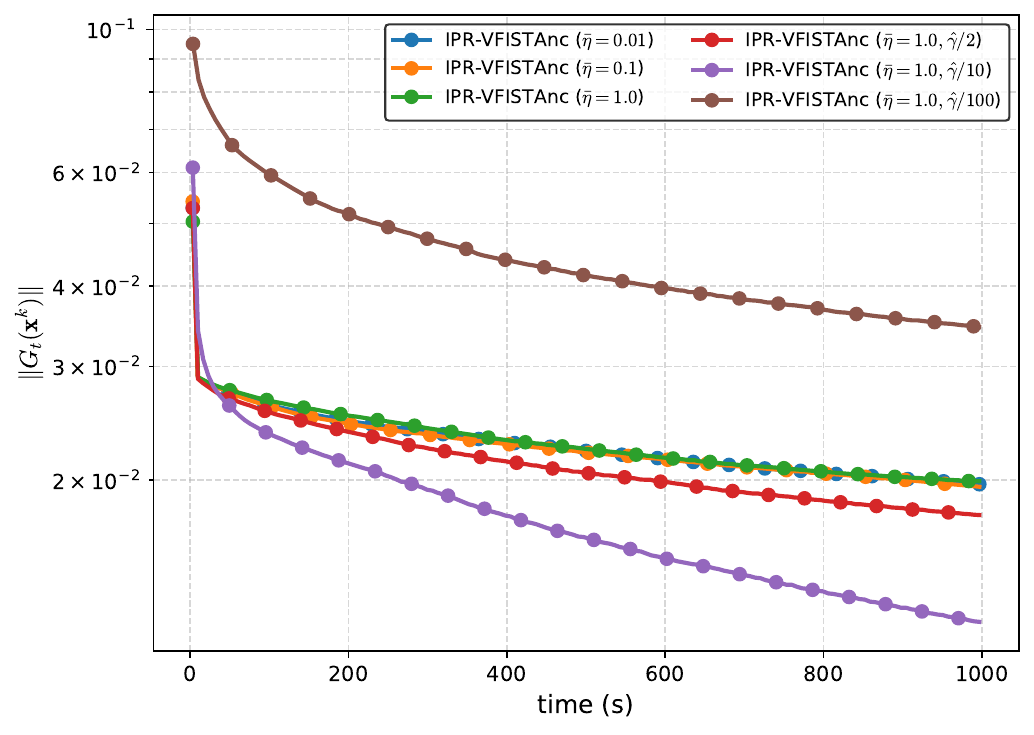}
            \caption{$\text{IPR-VFISTA}_{\text{nc}}$ configurations.}
            \label{fig:app_least_squares_ipr}
        \end{subfigure}
        \caption{Outer gradient mapping norm $\norm{G_{t}\left(\bx^{k}\right)}$ comparisons across tested hyperparameter configurations on the non-convex regularized least squares problem (averaged over 5 runs).}
        \label{fig:app_least_squares}
    \end{figure}

\begin{figure}[!htbp]
        \centering
        \begin{subfigure}[t]{0.48\textwidth}
            \centering
            \includegraphics[width=\textwidth]{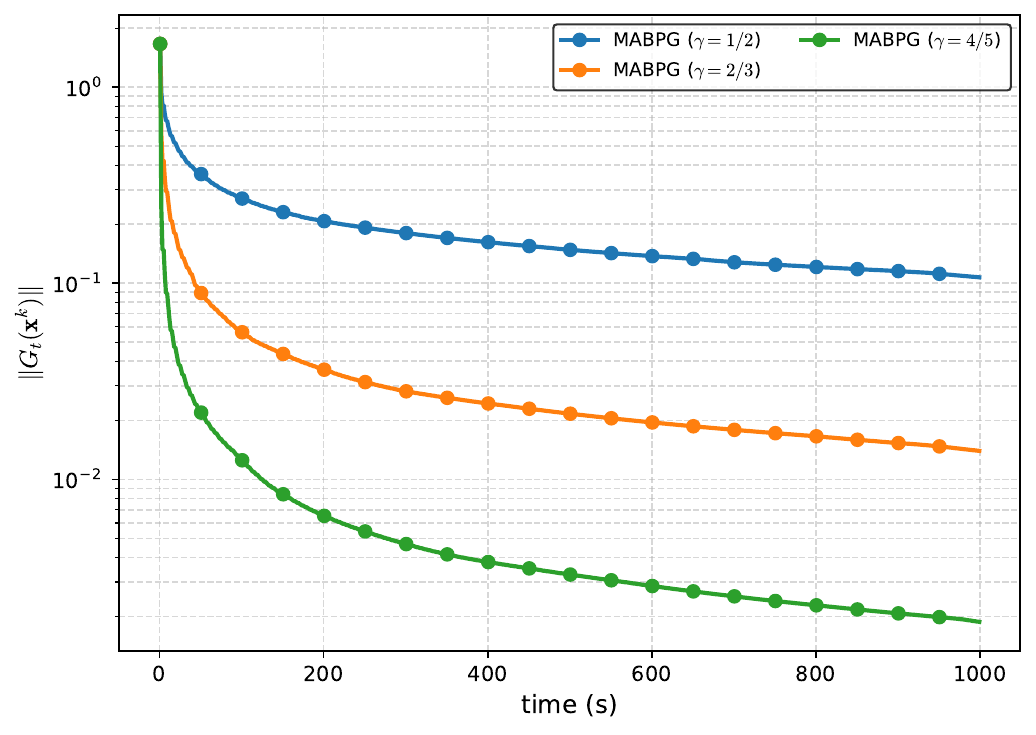}
            \caption{MABPG with $\gamma \in \{1/2, 2/3, 4/5\}$.}
            \label{fig:gamma_comparison_rcv1}
        \end{subfigure}
        \hfill
        \begin{subfigure}[t]{0.48\textwidth}
            \centering
            \includegraphics[width=\textwidth]{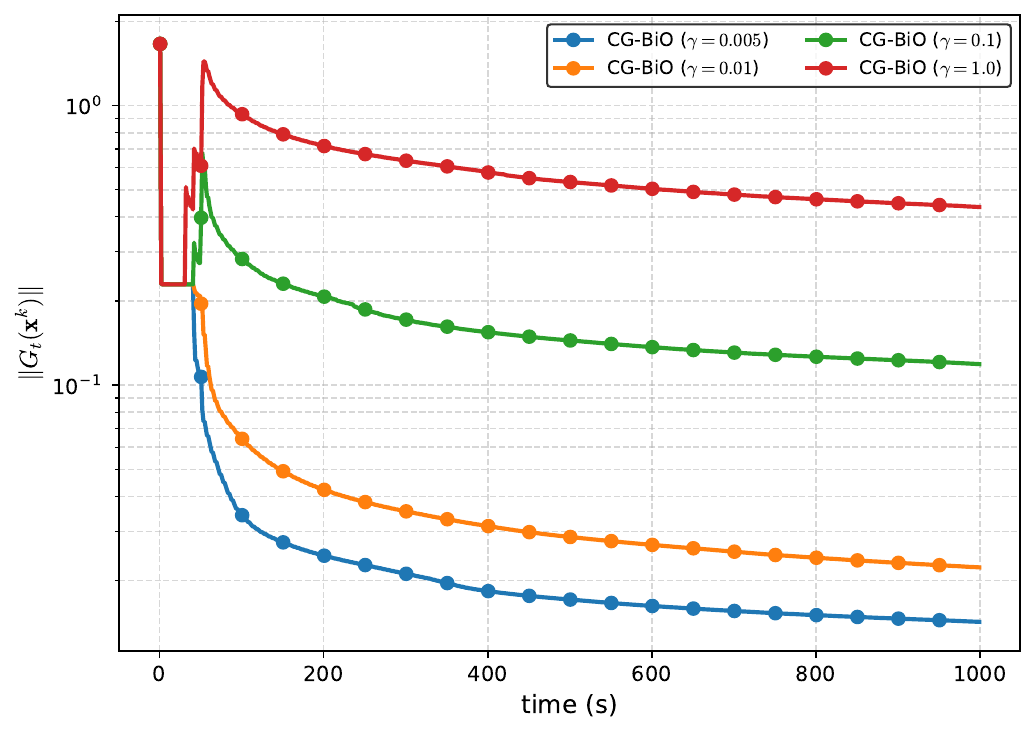}
            \caption{CG-BiO step-size coefficients.}
            \label{fig:app_rcv1_cgbio}
        \end{subfigure}
        \\[0.8em]
        \begin{subfigure}[t]{0.48\textwidth}
            \centering
            \includegraphics[width=\textwidth]{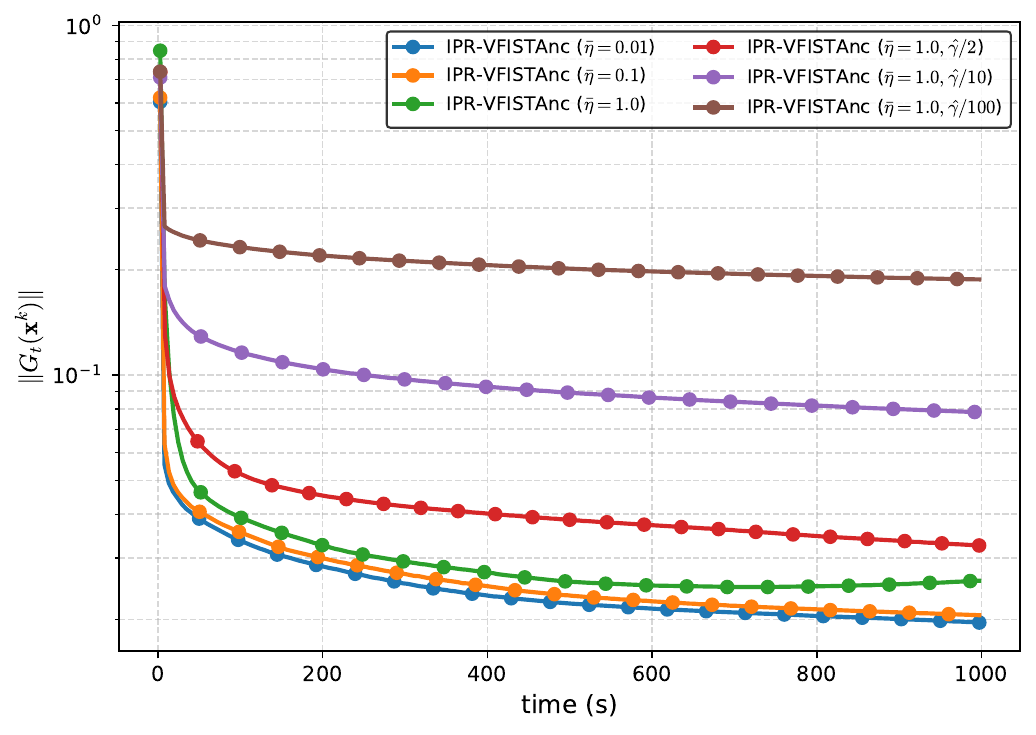}
            \caption{$\text{IPR-VFISTA}_{\text{nc}}$ configurations.}
            \label{fig:app_rcv1_ipr}
        \end{subfigure}
        \caption{Outer gradient mapping norm $\norm{G_{t}\left(\bx^{k}\right)}$ comparisons across tested hyperparameter configurations on the over-parameterized RCV1-v2 text classification problem (averaged over 5 runs).}
        \label{fig:app_rcv1}
    \end{figure}

    Figures \ref{fig:app_least_squares} and \ref{fig:app_rcv1} display $\norm{G_{t}\left(\bx^{k}\right)}$ versus elapsed time on the least-squares problem (Section \ref{Subsec:LeastSquares}) and the RCV1-v2 text classification problem (Section \ref{Subsec:TextClass}), respectively. For CG-BiO, $c = 0.005$ consistently achieves the lowest $\norm{G_{t}\left(\bx^{k}\right)}$ in both experiments. For vanilla MABPG, $\gamma = 4/5$ performs best on both problems, and is therefore inherited by MABPG-P. In contrast, $\text{IPR-VFISTA}_{\text{nc}}$ favors different settings across the two problems: $\bar{\eta} = 1$ with outer step size $h/10$ performs best on least squares, whereas $\bar{\eta} = 0.01$ with outer step size $h$ performs best on RCV1-v2.

\FloatBarrier
\section{Proofs of Propositions \ref{P:GradMapFWgap} and \ref{P:RegStepFW}} \label{App:StationarityBounds}
    \begin{proof}[Proof of Proposition \ref{P:GradMapFWgap}]
        Let $\by^{\ast} = \argmin_{\by \in X^{\ast}} \left\{ \act{\nabla\sigma\left(\bx\right) , \by - \bx} + \psi\left(\by\right) + \frac{1}{2t}\norm{\by - \bx}^{2} \right\}$. If $\by^{\ast} = \bx$, then $G_{t}^{\sigma , \psi + \delta_{X^{\ast}}}\left(\bx\right) = \bo$ (see \eqref{D:GradientMapping}) and \eqref{P:GradMapFWgap:1} holds trivially since $\pos{\GGG_{OP}\left(\bx\right)} \geq 0$. Assume henceforth that $\by^{\ast} \neq \bx$. Since $\bx^{\ast} \in X^{\ast}$, we have
        \begin{align}
            \act{\nabla \sigma\left(\bx\right) , \by^{\ast} - \bx} + \psi\left(\by^{\ast}\right) + \frac{1}{2t}\norm{\by^{\ast} - \bx}^{2} & \leq \act{\nabla \sigma\left(\bx\right) , \bx^{\ast} - \bx} + \psi\left(\bx^{\ast}\right) + \frac{1}{2t}\norm{\bx^{\ast} - \bx}^{2} \nonumber \\
            & \leq \act{\nabla \sigma\left(\bx\right) , \bx^{\ast} - \bx} + \psi\left(\bx\right) + \act{\xi , \bx^{\ast} - \bx} + \frac{1}{2t}\norm{\bx^{\ast} - \bx}^{2} \nonumber \\
            & \leq \norm{\nabla \sigma\left(\bx\right) + \xi} \cdot \norm{\bx^{\ast} - \bx} + \psi\left(\bx\right) + \frac{1}{2t}\norm{\bx^{\ast} - \bx}^{2}, \label{P:GradMapFWgap:2}
        \end{align}
        where the second inequality follows from the subgradient inequality for the convex function $\psi$ (recall that $\xi \in \partial \psi\left(\bx^{\ast}\right)$), and the last inequality follows from the Cauchy-Schwarz inequality. By the definition of the normalized FW gap (see Definition \ref{D:NormFWGap}), we have
        \begin{equation*}
            \act{\nabla \sigma\left(\bx\right) , \bx - \by^{\ast}} + \psi\left(\bx\right) - \psi\left(\by^{\ast}\right) \leq \GGG_{OP}\left(\bx\right) \cdot \norm{\by^{\ast} - \bx},
        \end{equation*}
        which combined with \eqref{P:GradMapFWgap:2} yields
        \begin{equation*}
            \frac{1}{2t}\norm{\by^{\ast} - \bx}^{2} \leq \GGG_{OP}\left(\bx\right) \cdot \norm{\by^{\ast} - \bx} + \norm{\nabla \sigma\left(\bx\right) + \xi} \cdot \norm{\bx^{\ast} - \bx} + \frac{1}{2t}\norm{\bx^{\ast} - \bx}^{2}.
        \end{equation*}
        Hence, we can rewrite the inequality as follows (after multiplying both sides by $2t$)
        \begin{equation*}
            \left(\norm{\by^{\ast} - \bx} - t\GGG_{OP}\left(\bx\right)\right)^{2} \leq t^{2}\GGG_{OP}\left(\bx\right)^{2} + 2t\norm{\nabla \sigma\left(\bx\right) + \xi} \cdot \norm{\bx^{\ast} - \bx} + \norm{\bx^{\ast} - \bx}^{2}.
        \end{equation*}
        Taking the square root of both sides, using $u \leq \left|u\right|$ on the left-hand side and applying twice the fact that $\sqrt{a + b} \leq \sqrt{a} + \sqrt{b}$ for all $a , b \geq 0$ on the right-hand side (with $\sqrt{t^{2}\GGG_{OP}\left(\bx\right)^{2}} = t\left|\GGG_{OP}\left(\bx\right)\right|$), we obtain
        \begin{equation*}
            \norm{\by^{\ast} - \bx} - t\GGG_{OP}\left(\bx\right) \leq t\left|\GGG_{OP}\left(\bx\right)\right| + \norm{\bx^{\ast} - \bx} + \sqrt{2t\norm{\nabla \sigma\left(\bx\right) + \xi} \cdot \norm{\bx^{\ast} - \bx}}.
        \end{equation*}
        Since $\GGG_{OP}\left(\bx\right) + \left|\GGG_{OP}\left(\bx\right)\right| = 2\pos{\GGG_{OP}\left(\bx\right)}$, rearranging yields
        \begin{equation*}
            \norm{\by^{\ast} - \bx} \leq 2t\pos{\GGG_{OP}\left(\bx\right)} + \norm{\bx^{\ast} - \bx} + \sqrt{2t\norm{\nabla \sigma\left(\bx\right) + \xi} \cdot \norm{\bx^{\ast} - \bx}}.
        \end{equation*}
        Since $\bx - \by^{\ast} = tG_{t}^{\sigma , \psi + \delta_{X^{\ast}}}\left(\bx\right)$ (see \eqref{D:GradientMapping}), dividing both sides by $t$ yields the desired result.
    \end{proof}

    \begin{proof}[Proof of Proposition \ref{P:RegStepFW}]
        From Proposition \ref{P:ProxGradFundProperty}(i), we have $G_{1/\beta}\left(\bv\right) - \nabla f_{\eta}\left(\bv\right) \in \partial g_{\eta}\left(\bx^{+}\right)$. Since $g_{\eta}$ is convex, for any $\bx \in \real^{n}$, it follows from the subgradient inequality that
        \begin{align*}
            g_{\eta}\left(\bx^{+}\right) - g_{\eta}\left(\bx\right) & \leq \act{G_{1/\beta}\left(\bv\right) - \nabla f_{\eta}\left(\bv\right) , \bx^{+} - \bx} \\
            & = \act{G_{1/\beta}\left(\bv\right) + \nabla f_{\eta}\left(\bx^{+}\right) - \nabla f_{\eta}\left(\bv\right) - \nabla f_{\eta}\left(\bx^{+}\right) , \bx^{+} - \bx} \\
            & \leq \norm{G_{1/\beta}\left(\bv\right)} \cdot \norm{\bx^{+} - \bx} + \beta \norm{\bx^{+} - \bv} \cdot \norm{\bx^{+} - \bx} - \act{\nabla f_{\eta}\left(\bx^{+}\right) , \bx^{+} - \bx} \\
            & = 2\norm{G_{1/\beta}\left(\bv\right)} \cdot \norm{\bx^{+} - \bx} - \act{\nabla f_{\eta}\left(\bx^{+}\right) , \bx^{+} - \bx} \\
            & \leq 2\norm{G_{1/\beta}\left(\bv\right)} \cdot \norm{\bx^{+} - \bx} + f\left(\bx\right) - f\left(\bx^{+}\right) - \eta\act{\nabla \sigma\left(\bx^{+}\right) , \bx^{+} - \bx},
        \end{align*}
        where the second inequality follows from the Cauchy-Schwarz inequality and the $\beta$-Lipschitz continuity of $\nabla f_{\eta}$, the last equality follows from the fact that $\bx^{+} = \bv - \beta^{-1}G_{1/\beta}\left(\bv\right)$, and the last inequality follows from the convexity of $f$ and the definition of $f_{\eta}$. Thus, rearranging the obtained inequality and using the definition of $g_{\eta}$, yields \eqref{P:RegStepFW:1}.
\medskip

        Moreover, for any $\bx \in X^{\ast}$, we have that $\varphi\left(\bx^{+}\right) - \varphi\left(\bx\right) \geq 0$, and thus
        \begin{equation*}
            \eta\left(\psi\left(\bx^{+}\right) - \psi\left(\bx\right) + \act{\nabla \sigma\left(\bx^{+}\right) , \bx^{+} - \bx}\right) \leq 2\norm{G_{1/\beta}\left(\bv\right)} \cdot \norm{\bx^{+} - \bx}.
        \end{equation*}
        Hence, for any $\bx \in X^{\ast} \setminus \{ \bx^{+} \}$, dividing both sides in $\eta\norm{\bx^{+} - \bx}$ yields
        \begin{equation*}
            \frac{\psi\left(\bx^{+}\right) - \psi\left(\bx\right) + \act{\nabla \sigma\left(\bx^{+}\right) , \bx^{+} - \bx}}{\norm{\bx^{+} - \bx}} \leq \frac{2}{\eta} \cdot \norm{G_{1/\beta}\left(\bv\right)}.
        \end{equation*}
        Therefore, the result now follows by taking the supremum over all $\bx \in X^{\ast} \setminus \{ \bx^{+} \}$.
    \end{proof} 

\section{Proof of Theorem \ref{T:TotalRate}} \label{App:ProofTotalRate}
    Throughout this appendix, $\bz = \bx^{k_{best} + 1}$ is the output of Phase 1 of MABPG-P, and $\bar{\bz} = P_{X^{\ast}}\left(\bz\right)$ is its projection onto $X^{\ast}$, that is, the unique solution of Problem \eqref{P:ProjectionBiLevelProblem}. We recall the guarantees of FBi-PG for Problem \eqref{P:ProjectionBiLevelProblem}, which follow from \cite[Theorem 4]{MST2026} with the outer objective function $h_{\bz}$, since $\varphi$ satisfies the H\"olderian error bound with parameter $\tau$ (see Assumption \ref{A:AssumpB}(i)). For all $s \geq 1$, we have
    \begin{equation} \label{App:FBiPG}
        \varphi\left(\bw^{s}\right) - \varphi\left(\bar{\bz}\right) \leq \frac{aC_{F}}{\left(s + 1\right)^{2}} \quad \text{and} \quad \norm{\bw^{s} - \bz}^{2} - \norm{\bar{\bz} - \bz}^{2} \leq \frac{2C_{F}}{\left(s + 1\right)^{1/2}},
    \end{equation}
    where $C_{F} := a^{2}\norm{\bar{\bz} - \bz}^{2}\max\left\{ \tilde{\beta} + 1/2 , \, 4a^{3}/\tau \right\}$.
\medskip

    \noindent We begin with the effect of an approximate projection on the gradient mapping.
    \begin{lemma} \label{L:BoundGradientMappingProjection}
        Let $t > 0$, $\varepsilon > 0$ and $\bx \in \real^{n}$. Then, for any $\by \in \real^{n}$ with $\norm{\by - P_{X^{\ast}}\left(\bx\right)} \leq \varepsilon$, we have
        \begin{equation*}
            \norm{G_{t}^{\sigma , \psi + \delta_{X^{\ast}}}\left(\by\right)} \leq \left(3 + t\beta_{\sigma}\right)\norm{G_{t}^{\sigma , \psi + \delta_{X^{\ast}}}\left(\bx\right)} + \frac{2 + t\beta_{\sigma}}{t}\varepsilon.
        \end{equation*}
    \end{lemma}
    \begin{proof}
        For simplicity, we use $G_{t} := G_{t}^{\sigma , \psi + \delta_{X^{\ast}}}$ and $T_{t} := T_{t}^{\sigma , \psi + \delta_{X^{\ast}}}$. Since the proximal mapping is nonexpansive \cite[Theorem 6.42]{B2017-B} and $\nabla\sigma$ is $\beta_{\sigma}$-Lipschitz continuous, for any $\bu , \bv \in \real^{n}$ we have (recall that $T_{t}\left(\bu\right) = \prox_{t\left(\psi + \delta_{X^{\ast}}\right)}\left(\bu - t\nabla\sigma\left(\bu\right)\right)$ for all $\bu \in \real^{n}$)
        \begin{equation*}
            \norm{T_{t}\left(\bu\right) - T_{t}\left(\bv\right)} \leq \norm{\bu - \bv} + t\norm{\nabla\sigma\left(\bu\right) - \nabla\sigma\left(\bv\right)} \leq \left(1 + t\beta_{\sigma}\right)\norm{\bu - \bv}.
        \end{equation*}
        Therefore, from the definition of the gradient mapping (see \eqref{D:GradientMapping}), we obtain
        \begin{equation} \label{L:BoundGradientMappingProjection:1}
            \norm{G_{t}\left(\bu\right) - G_{t}\left(\bv\right)} \leq \frac{1}{t}\left(\norm{\bu - \bv} + \norm{T_{t}\left(\bu\right) - T_{t}\left(\bv\right)}\right) \leq \frac{2 + t\beta_{\sigma}}{t}\norm{\bu - \bv}.
        \end{equation}
        Since $T_{t}\left(\bx\right) \in X^{\ast}$, the definition of the projection yields $\norm{\bx - P_{X^{\ast}}\left(\bx\right)} \leq \norm{\bx - T_{t}\left(\bx\right)} = t\norm{G_{t}\left(\bx\right)}$, and thus
        \begin{equation*}
            \norm{\by - \bx} \leq \norm{\by - P_{X^{\ast}}\left(\bx\right)} + \norm{P_{X^{\ast}}\left(\bx\right) - \bx} \leq \varepsilon + t\norm{G_{t}\left(\bx\right)}.
        \end{equation*}
        Using \eqref{L:BoundGradientMappingProjection:1} with $\bu = \by$ and $\bv = \bx$, we get
        \begin{equation*}
            \norm{G_{t}\left(\by\right)} \leq \norm{G_{t}\left(\bx\right)} + \frac{2 + t\beta_{\sigma}}{t}\left(\varepsilon + t\norm{G_{t}\left(\bx\right)}\right) = \left(3 + t\beta_{\sigma}\right)\norm{G_{t}\left(\bx\right)} + \frac{2 + t\beta_{\sigma}}{t}\varepsilon.
        \end{equation*}
        This completes the proof.
    \end{proof}
    The next lemma shows that Phase 2 brings $\bw^{s}$ close to $\bar{\bz}$.
    \begin{lemma} \label{L:Phase2DistanceBound}
        For all $s \in \nn$, we have
        \begin{equation*}
            \norm{\bw^{s} - \bar{\bz}}^{2} \leq \frac{2C_{F}}{\left(s + 1\right)^{1/2}} + \frac{2}{s + 1}\sqrt{\frac{aC_{F}}{\tau}} \cdot \norm{\bz - \bar{\bz}}.
        \end{equation*}
    \end{lemma}
    \begin{proof}
        Expanding $\norm{\bw^{s} - \bz}^{2} = \norm{\left(\bw^{s} - \bar{\bz}\right) + \left(\bar{\bz} - \bz\right)}^{2}$, we obtain
        \begin{equation} \label{L:Phase2DistanceBound:1}
            \norm{\bw^{s} - \bar{\bz}}^{2} = \norm{\bw^{s} - \bz}^{2} - \norm{\bar{\bz} - \bz}^{2} + 2\act{\bw^{s} - \bar{\bz} , \bz - \bar{\bz}}.
        \end{equation}
        From \eqref{App:FBiPG}, the first term is bounded by $2C_{F}/\left(s + 1\right)^{1/2}$. For the last term, since $P_{X^{\ast}}\left(\bw^{s}\right) \in X^{\ast}$ and $\bar{\bz} = P_{X^{\ast}}\left(\bz\right)$, the projection theorem yields $\act{P_{X^{\ast}}\left(\bw^{s}\right) - \bar{\bz} , \bz - \bar{\bz}} \leq 0$, and thus, using the Cauchy-Schwarz inequality, we have
        \begin{equation*}
            \act{\bw^{s} - \bar{\bz} , \bz - \bar{\bz}} \leq \act{\bw^{s} - P_{X^{\ast}}\left(\bw^{s}\right) , \bz - \bar{\bz}} \leq \dist\left(\bw^{s} , X^{\ast}\right) \cdot \norm{\bz - \bar{\bz}}.
        \end{equation*}
        From Assumption \ref{A:AssumpB}(i) and \eqref{App:FBiPG}, we have
        \begin{equation*}
            \tau \cdot \dist\left(\bw^{s} , X^{\ast}\right)^{2} \leq \varphi\left(\bw^{s}\right) - \varphi\left(\bar{\bz}\right) \leq \frac{aC_{F}}{\left(s + 1\right)^{2}}.
        \end{equation*}
        Plugging these bounds into \eqref{L:Phase2DistanceBound:1} completes the proof.
    \end{proof}
    We are now ready to prove Theorem \ref{T:TotalRate}.
    \begin{proof}[Proof of Theorem \ref{T:TotalRate}]
        From Theorem \ref{T:ImprovedInnerRate} with $\gamma = 1/2$, we have $\norm{\bz - \bar{\bz}} \leq \left(2C + W\right)/\left(\tau k^{1/2}\right)$, and thus
        \begin{equation} \label{App:ProofTotalRate:1}
            C_{F} \leq \frac{H}{k}, \quad \text{where} \quad H := \frac{a^{2}\left(2C + W\right)^{2}}{\tau^{2}}\max\left\{ \tilde{\beta} + \frac{1}{2} , \, \frac{4a^{3}}{\tau} \right\}.
        \end{equation}
        Since $\bar{\bz} \in X^{\ast}$, we have $\varphi\left(\bar{\bz}\right) = \varphi\left(\bx\right)$. Therefore, from \eqref{App:FBiPG} with $s = 2k$ and \eqref{App:ProofTotalRate:1}, we obtain
        \begin{equation*}
            \varphi\left(\bw^{2k}\right) - \varphi\left(\bx\right) \leq \frac{aC_{F}}{\left(2k + 1\right)^{2}} \leq \frac{aH}{4k^{3}},
        \end{equation*}
        which proves the first bound. From Lemma \ref{L:Phase2DistanceBound} with $s = 2k$, \eqref{App:ProofTotalRate:1}, $2k + 1 \geq 2k$ and $\norm{\bz - \bar{\bz}} \leq \left(2C + W\right)/\left(\tau k^{1/2}\right)$, we obtain
        \begin{equation*}
            \norm{\bw^{2k} - \bar{\bz}}^{2} \leq \frac{2H}{k\left(2k\right)^{1/2}} + \frac{1}{k}\sqrt{\frac{aH}{\tau k}} \cdot \frac{2C + W}{\tau k^{1/2}} = \frac{\sqrt{2}H}{k^{3/2}} + \frac{2C + W}{\tau}\sqrt{\frac{aH}{\tau}} \cdot \frac{1}{k^{2}} \leq \frac{H'}{k^{3/2}},
        \end{equation*}
        where $H' := \sqrt{2}H + \frac{2C + W}{\tau}\sqrt{aH/\tau}$, and the last inequality follows from the fact that $k^{-2} \leq k^{-3/2}$ for $k \geq 1$. Now, from Lemma \ref{L:BoundGradientMappingProjection} with $t = 1/\beta_{\sigma}$, $\bx = \bz$, $\by = \bw^{2k}$ and $\varepsilon = \sqrt{H'}k^{-3/4}$, and from Theorem \ref{T:GradientMappingRate} with $\gamma = 1/2$, we obtain
        \begin{equation*}
            \norm{G_{1/\beta_{\sigma}}^{\sigma , \psi + \delta_{X^{\ast}}}\left(\bw^{2k}\right)} \leq \frac{4D}{k^{1/4}} + 3\beta_{\sigma} \cdot \frac{\sqrt{H'}}{k^{3/4}} \leq \frac{\Theta_{2}}{k^{1/4}},
        \end{equation*}
        where the last inequality follows from the fact that $k^{-3/4} \leq k^{-1/4}$ for $k \geq 1$. This completes the proof with
        \begin{equation} \label{App:ProofTotalRate:Theta}
            \Theta_{1} := \frac{aH}{4} \quad \text{and} \quad \Theta_{2} := 4D + 3\beta_{\sigma}\sqrt{H'},
        \end{equation}
        where $D$ is defined in Theorem \ref{T:GradientMappingRate}, and $H$ and $H'$ are defined above.
    \end{proof}

\bibliographystyle{plain}
\bibliography{notes}

\end{document}